\documentclass[11pt]{article}
\usepackage{amscd, amsmath, amssymb}

\title{Cheeger manifolds and the classification of biquotients}
\author{Burt Totaro}
\date{  }

\def\Z{\text{\bf Z}}
\def\Q{\text{\bf Q}}
\def\R{\text{\bf R}}
\def\C{\text{\bf C}}
\def\P{\text{\bf P}}
\def\F{\text{\bf F}}
\def\arrow{\rightarrow}
\def\inj{\hookrightarrow}

\def\qed{\ QED \vspace{0.2in}}

\def\HH{\text{\bf H}}
\def\Ca{\text{\bf Ca}}
\def\back{\negthinspace \setminus}
\def\quot{/\negthinspace /}

\def\odd{\text{odd}}
\def\til{\widetilde}

\begin{document}
\maketitle

\newtheorem{theorem}{Theorem}[section]
\newtheorem{corollary}[theorem]{Corollary}
\newtheorem{lemma}[theorem]{Lemma}
\newtheorem{tab}[theorem]{Table}
\newtheorem{convention}[theorem]{Convention}

A closed manifold is called a biquotient if it is diffeomorphic
to $K\back G/H$ for some compact Lie group $G$ with closed subgroups
$K$ and $H$ such that $K$ acts freely on $G/H$. Every biquotient has
a Riemannian metric of nonnegative sectional curvature. In fact,
almost all known manifolds of nonnegative curvature are biquotients.
The only known closed manifolds of nonnegative sectional curvature which were
not defined in this way
are those found by Cheeger \cite{Cheeger} in 1973,
and Grove and Ziller \cite{GZ} in 2000.

In order to understand these constructions better, this paper analyzes
which of the Cheeger and Grove-Ziller
manifolds are actually diffeomorphic to biquotients. In the process,
we develop a general approach to the classification of biquotient
manifolds. A priori, it is hard
to determine whether a given manifold is a biquotient, because
there is no obvious upper bound
on the dimension of the groups involved. We give a procedure which
allows us to reduce the dimension of the groups needed to describe
a given manifold as a biquotient.
A surprising application is
that Gromoll-Meyer's example of an exotic 7-sphere which is a biquotient
\cite{GM}
is the only exotic sphere of any dimension which is a biquotient.
More generally, we classify
all biquotients which are simply connected rational homology spheres
of any dimension (Theorem \ref{sphere}). The classification of
biquotients which are simply connected rational homology spheres has been
given independently by Kapovitch and
Ziller \cite{KZ}. More generally,
they classify all simply connected biquotients whose rational cohomology
ring is generated by one element.

Another application of our general classification theory for biquotients
is the precise determination of which Cheeger manifolds,
the connected sums
of two rank-one symmetric spaces with any orientations,
are diffeomorphic
to biquotients (Theorem \ref{cheeger}).
 Some of the Cheeger manifolds were known to be biquotients,
such as $\C\P^2 \# -\C\P^2$ (the nontrivial $S^2$-bundle over $S^2$),
but we find that many of the other
Cheeger manifolds are also biquotients, such as $\C\P^2\# \C\P^2$. On the other
hand, the Cheeger manifold $\C\P^4 \# \HH\P^2$ is not diffeomorphic
(or even homotopy equivalent) to a biquotient. The positive result
that some Cheeger manifolds such as $\C\P^2\# \C\P^2$ are biquotients
implies that they have nonnegatively curved Riemannian metrics
with various good properties that were not clear from Cheeger's
construction.
First, the new metrics are real analytic. Further, the new
metrics determine a natural complex-analytic
structure on the whole tangent bundle, since Aguilar \cite{Aguilar}
showed that biquotients have this property;
my paper \cite{Totarocomplex} has a weaker result which applies
to all the Cheeger manifolds.
Finally, the geodesic flow for the new metrics on the Cheeger manifolds
has zero topological entropy. Paternain \cite{Paternain} proved
the latter property for Cheeger's metric on $\C\P^2\# \C\P^2$, but
now we know it for some metric on a large class of the Cheeger manifolds.

Finally, our general classification results imply that there are only
finitely many diffeomorphism classes of 2-connected biquotients
in any given dimension (Theorem \ref{finite}).
This fails for biquotients that are only simply
connected
already in dimension 6, as I showed in \cite{Totarocurv}.
Also, the finiteness of 2-connected biquotients shows the distance
between biquotients and general manifolds of nonnegative sectional
curvature, since Grove and Ziller have  constructed
nonnegatively curved metrics on all $S^3$-bundles over $S^4$,
giving infinitely many homotopy types of 2-connected 7-manifolds
with nonnegative sectional curvature \cite{GZ}.

I would like to thank Gabriel Paternain for many useful conversations.
Also, thanks to Wolfgang Ziller for several references to earlier work.

\tableofcontents

\section{Notation}

We begin with the equivalence of several definitions of biquotient
manifolds, pointed out by Eschenburg (\cite{Eschenburgexample},
\cite{Eschenburgcoho}).

\begin{lemma}
The following properties of a closed smooth manifold $M$ are equivalent.

(1) $M=K\back G/H$ for some compact Lie group $G$ with closed subgroups
$K$ and $H$ such that $K$ acts freely on $G/H$.

(2) $M=G/H$ for some compact Lie groups $G$ and $H$ together with a
homomorphism $H\arrow G\times G$ such that $H$ acts freely on $G$
by left and right translation, $(g_1,g_2)(g):=g_1gg_2^{-1}$.

(3) $M=G/H$ for some compact Lie groups $G$ and $H$ together with
a homomorphism $H\arrow (G\times G)/Z(G)$ such that $H$ acts
freely on $G$.
\end{lemma}

In (3) we are using that the center $Z(G)$, imbedded diagonally in $G\times
G$,
acts trivially on $G$ by  left and right translation.

{\bf Proof. }Clearly (1) implies (2) and (2) implies (3). We prove that (3)
implies (1). The point is that $(G\times G)/Z(G)$ acts transitively on $G$
with stabilizer at $1\in G$ equal to the diagonal subgroup $G/Z(G)$. So,
if $M=G/H$ as in (3), then we can also describe $M$ as in (1) by
$$M=\big( G/Z(G)\big) \back \big( (G\times G)/Z(G)\big) /H.$$
\qed

In the paper, we use definition (3) of biquotients; that is, ``a biquotient
$G/H$''
will mean that $G$ and $H$ are compact Lie groups and
we are given a homomorphism $H\arrow (G\times G)/Z(G)$
such that $H$ acts freely on $G$. When $G$ and $H$ are simply connected,
such a homomorphism lifts uniquely to $G\times G$.

We call
a connected compact Lie group {\it simple }if it is nonabelian and every
proper
normal subgroup is finite.
We write $Sp(2a)$ for the simply connected simple group of type $C_a$,
which topologists usually call $Sp(a)$.
We define a
{\it simple factor }of a connected compact Lie group $H$ to be the universal
cover of a simple normal subgroup of $H$.

By definition, the {\it Dynkin index }of a homomorphism $H\arrow G$
of simply connected simple groups is the integer corresponding to the
homomorphism $\pi_3H\arrow \pi_3G$, both groups being canonically
isomorphic to $\Z$. Dynkin computed the Dynkin index in many cases
(\cite{Dynkin}, Chapter I, section 2).
Finally, we write $UT(S^n)$ for the unit tangent bundle
of the $n$-sphere, $UT(S^n)=Spin(n+1)/Spin(n-1)$.

\section[Cheeger manifolds that are biquotients]{The Cheeger manifolds
that are diffeomorphic to biquotients}
\label{positive}

By definition, a Cheeger manifold is the connected sum of any two
rank-one
symmetric spaces with any orientations. The rank-one symmetric spaces,
besides the sphere which is not interesting for this purpose, are the real,
complex, and quaternionic projective spaces, together with the Cayley plane
associated to the octonions.
The Cheeger manifolds are the only known examples of connected sums,
with neither summand a homotopy sphere and at least one summand
not a rational homology sphere,
which admit metrics of nonnegative sectional curvature. In fact,
the conjecture that manifolds of nonnegative sectional curvature are
elliptic in the sense of F\'elix, Halperin, and Thomas \cite{FHTell}
would imply that any connected sum which admits
a metric of nonnegative sectional curvature must be roughly of the Cheeger
type. Precisely,
a connected sum $M_1\# M_2$ of simply connected manifolds
which is elliptic, and such that
neither $M_1$ nor $M_2$ is a $k$-homology sphere for a field $k$,
must have the rings
$H^*(M_1,k)$ and $H^*(M_2,k)$ both generated by a single element,
as follows from Lambrechts \cite{Lambrechts}, Theorem 3. Combining
this with Adams and Atiyah's results on the Hopf invariant problem
\cite{AA} shows that if a connected sum of simply connected
manifolds is elliptic, with neither summand a homotopy sphere
and at least one summand not a rational homology sphere,
then both summands have the integral cohomology
ring of $\C\P^n$, $\HH\P^n$, or $\Ca\P^2$.

\begin{theorem}
\label{cheeger}
The following Cheeger manifolds are diffeomorphic to biquotients.
First, $\C\P^n\# -\C\P^n$,
$\HH\P^n\# -\HH\P^n$, and $\Ca\P^2\# -\Ca\P^2$.
Next, $\R\P^n\# \R\P^n$, $\R\P^{2n}\# \C\P^n$, $\R\P^{4n}\# \HH\P^n$,
and $\R\P^{16}\# \Ca\P^2$. Here $\R\P^n$ is non-orientable for $n$ odd and
has an orientation-reversing diffeomorphism for $n$ even, so the
orientations
of the summands do not matter in these cases. Next, $\C\P^{2n}\# -\HH\P^n$,
$\HH\P^4\# -\Ca\P^2$, and $\C\P^8\# -\Ca\P^2$.
Finally, $\C\P^n\# \C\P^n$, $\HH\P^n\# \HH\P^n$, and
$\C\P^{4e+2}\# \HH\P^{2e+1}$.

The remaining Cheeger manifolds are not even homotopy equivalent to
biquotients.
These are $\Ca\P^2\# \Ca\P^2$,
$\C\P^8\# \Ca\P^2$, $\HH\P^4\# \Ca\P^2$,
and $\C\P^{4e}\# \HH\P^{2e}$.
\end{theorem}

{\bf Proof. }In this section, we prove only the first, positive, statement.
We will
prove the negative statement in sections \ref{16} and \ref{negative},
after setting up
a general classification theory of biquotients.

We begin with the cases which are
straightforward generalizations of Cheeger's observation
that $\C\P^2\# -\C\P^2$ is a biquotient \cite{Cheeger}. Let $A_k$ denote
the standard algebra ($\R$, $\C$, $\HH$, or the octonions $\Ca$) of
dimension
$k$
over $\R$, where $k=1$, 2, 4, or 8. Then, for $k=1$, 2, or 4,
$A_k\P^n\# -A_k\P^n$ is an $S^k$-bundle over $A_k\P^{n-1}$, namely
the biquotient manifold $(S^{nk-1}\times S^k)/S^{k-1}$, where
$S^{k-1}$ is a group for $k=1$, 2, or 4, acting freely on $S^{nk-1}$
and by rotations on $S^k$. For $k=8$, $S^7$ is not a group, but we can
still describe $\Ca\P^2\# -\Ca\P^2$ as a biquotient,
$$(Spin(9)\times S^8)/Spin(8),$$
where $Spin(8)\arrow Spin(9)$ is the standard inclusion and $Spin(8)$
acts on $S^8\subset \R^9$ by the direct sum of an 8-dimensional
real spin representation and the trivial representation. Since
$Spin(9)/Spin(8)
=S^8$, this description exhibits $\Ca\P^2\# -\Ca\P^2$ as an $S^8$-bundle
over $S^8$.

Next, we check that the connected sum of $\R\P^n$ with any rank-one
symmetric space is diffeomorphic to a biquotient. As mentioned in the
theorem,
orientations do not matter in this case, because $\R\P^n$ is
non-orientable
for $n$ odd and has an orientation-reversing diffeomorphism for $n$ even.
For any closed $n$-manifold
$M$, the connected sum $\R\P^n\# M$ is doubly covered by $M\# -M$.
This suggests a way to view $\R\P^{nk}\# A_k\P^n$ as a biquotient: we
replace
$S^k$ in the above description of $A_k\P^n\# -A_k\P^n$ by $\R\P^k$.
That is:
\begin{align*}
\R\P^n\# \R\P^n &= (S^{n-1}\times S^1)/\Z/2 \\
\R\P^{2n}\# \C\P^n &= (S^{2n-1}\times \R\P^2)/S^1 \\
\R\P^{4n} \# \HH\P^n &= (S^{4n-1}\times \R\P^4)/S^3 \\
\R\P^{16} \# \Ca\P^2 &= (Spin(9)\times\R\P^8)/Spin(8)
\end{align*}

The manifolds
$\C\P^{2n}$, $\HH\P^{2n}$, and $\Ca\P^2$ have natural
orientations, corresponding
to the highest power of any generator of $H^2$, $H^4$, or $H^8$,
respectively.
In fact, $\HH\P^n$ has a natural orientation for all $n\geq 2$, but this
takes
more care to define. We use that the mod 3 Steenrod operation $P^1$ on
 $H^*(BSU(2),\F_3)$ acts by $P^1c_2=c_2^2$, as one checks by restricting
to the torus $S^1\subset SU(2)$. Since $H^4(\HH\P^n,\Z)$ has a generator
which is $c_2$ of an $SU(2)$-bundle over $\HH\P^n$, this generator $c_2$
satisfies
$P^1c_2=c_2^2$ in $H^8(\HH\P^n,\F_3)$.
 It follows that (for $n\geq 2$) there is a unique generator
$z$ of $H^4(\HH\P^n,\Z)$ such that $P^1z=-z^2$ in $H^8(\HH\P^n,\F_3)$,
namely $z=-c_2$. We define the natural orientation on $\HH\P^n$ to be the
one corresponding to the highest power $z^n$ of this generator $z$.

Using this orientation, we can define one more class of Cheeger manifolds
as bundles. The manifold $\C\P^{2n}\# -\HH\P^n$ is the $\C\P^2$-bundle
over $\HH\P^{n-1}$ associated to the rank-3 complex vector bundle
$E\oplus \C$, where $E$ is the tautological rank-2 complex vector bundle
over $\HH\P^{n-1}$. Thus $\C\P^{2n}\# -\HH\P^n$ is a biquotient of the form
$(\C\P^2\times S^{4n-1})/SU(2)$. We might try to imitate this construction
by viewing  $\HH\P^4\# -\Ca\P^2$ as an $\HH\P^2$-bundle over $S^8$ and
$\C\P^8\# -\Ca\P^2$ as a $\C\P^4$-bundle over $S^8$, but that turns out
to be impossible. For these constructions,
we would need a complex vector bundle $E$ over $S^8$ with
$c_4E$ equal to plus or minus the class of a point in $H^8(S^8,\Z)$,
whereas in fact every complex vector bundle on $S^8$ has $c_4E$ a multiple
of $(4-1)!=6$, by Bott periodicity.

Instead, we construct $\HH\P^4\# -\Ca\P^2$ as the quotient of an
$S^{11}$-bundle
over $S^8$ by a free $SU(2)$-action. Let $S^{-}:Spin(8)\arrow SO(8)$ denote
one of the two spin representations. We also write $S^{-}$ for the
associated
$Spin(9)$-equivariant real vector bundle of rank 8 over
$Spin(9)/Spin(8)=S^8$. Let $N$ be the $S^{11}$-bundle
over $S^8$ defined as the unit sphere
bundle
$S(S^{-}\oplus\R^4)$.
Define a homomorphism $SU(2)\arrow Spin(9)$ by $SU(2)\cong Spin(3)
\subset Spin(9)$. Using this homomorphism,
$SU(2)$ acts on $S^8$ and acts compatibly on the vector bundle $S^{-}$ over
$S^8$. Let $SU(2)$ act on $N$ by the given action on $S^8$ and on the vector
bundle
$S^{-}$, and by the standard faithful representation $V_{\R}$ of $SU(2)$ on
$\R^4$.

This action of $SU(2)$ on $N$ is free. To check this, it suffices to check
that
$SU(2)$ acts freely on the $S^7$-bundle $S(S^{-})$ over $S^8$ and on the
$S^3$-bundle $S(\R^4)=S^3\times S^8$ over $S^8$. The second statement is
clear by the choice of $SU(2)$-action on $\R^4$. To prove the first
statement,
first note that $Spin(9)$ acts on $S(S^{-})\cong S^{15}$ by the
spin representation of $Spin(9)$. Then use that the restriction of any
spin representation of $Spin(n)$ to $Spin(n-1)$ is a sum of spin
representations
of $Spin(n-1)$. Thus the action of $SU(2)\cong Spin(3)$ on $S(S^{-})=S^{15}$
is by a sum of copies of the 4-dimensional real spin representation
$V_{\R}$. It follows that $SU(2)$ acts freely on $S^{15}$.

Let $M$ be the quotient manifold $N/SU(2)$, of dimension 16. Clearly $M$
is a biquotient of the form $(Spin(9)\times S^{11})/(Spin(8)\times SU(2))$.
 By construction, $M$ is the union
of a disc bundle over $S(S^{-})/SU(2)=S^{15}/SU(2)=\HH\P^3$ and a disc
bundle
over $S(\R^4)/SU(2)=(S^8\times S^3)/SU(2)=S^8$ along their common boundary,
$S^{15}$.
So $M$ is  diffeomorphic to the connected sum of $\HH\P^4$ and $\Ca\P^2$,
with some orientations. To pin down the orientations, it is convenient to
compute
the cohomology ring of $M$. Since $N$ is an $SU(2)$-equivariant
$S^{11}$-bundle
over $S^8$, the quotient $M$ is the total space of a fibration
$S^{11}\arrow M\arrow S^8\quot SU(2)$ up to homotopy, where $S^8\quot SU(2)
=(S^8\times ESU(2))/SU(2)$ denotes the homotopy quotient via the given
homomorphism
$SU(2)\arrow Spin(9)$. 

We compute that
$$H^*(S^8\quot SU(2),\Z)=\Z[y,\chi(S^{-})]/((\chi(S^{-})+y^2)^2=0),$$
where $y$ in $H^4$ is $c_2$ of the standard representation $V$ of
$SU(2)$
and $\chi(S^{-})$ in $H^8$
is the Euler class of the $SU(2)$-equivariant
vector bundle $S^{-}$ over $S^8$, with some orientation. To check this,
note that the spectral sequence for the fibration $S^8\arrow S^8\quot SU(2)
\arrow BSU(2)$ degenerates, since all the cohomology is in even degrees.
This implies that the cohomology of $S^8\quot SU(2)$ has the above
form, for some relation of the form $\chi(S^{-})^2+b\chi(S^{-})y^2
+cy^4=0$ with $b,c\in\Z$. Furthermore, we have $TS^8\oplus \R\cong \R^9$
as $SU(2)$-equivariant vector bundles on $S^8$. It follows that
$$\chi(TS^8)^2=c_8(TS^8\otimes_{\R}\C)=0$$
in $H^{16}(S^8\quot SU(2),\Z)$. Also, $\chi(TS^8)$ has degree
2 on $S^8$ while $\chi(S^{-})$ has degree $\pm 1$ on $S^8$, so we must
have $\chi(TS^8)=\pm 2 \chi(S^{-})+dy^2$ in $H^8(S^8\quot SU(2),\Z)$
for some integer $d$. So $(\pm \chi(S^{-})+(d/2)y^2)^2=0$ in
$H^{16}(S^8\quot SU(2),\Q)$. By what we know about the form of the
relation, $d/2$ must be an integer $a$, and we have 
$$(\pm \chi(S^{-})+ay^2)^2 =0$$
in $H^{16}(S^8\quot SU(2),\Z)$. Finally, the action of $SU(2)
=Spin(3)\subset Spin(9)$ on $S^8$ preserves a 2-sphere, and so we have
an inclusion $BS^1\simeq S^2\quot SU(2)\arrow S^8\quot SU(2)$. 
Let $x$ be the standard generator of the polynomial ring
$H^*(BS^1,\Z)$ in degree 2. We compute that
the restriction map takes $y=c_2V$ to $-x^2$ and $\chi(S^{-})$
to $\pm x^4$, using that the restriction of the spin representation
$S^{-}$ to $S^1=Spin(2)\subset Spin(8)$ is the direct sum of 4 copies
of the standard 2-dimensional real representation of $S^1$.
Therefore, for a suitable orientation on $S^{-}$, the relation
in $H^{16}(S^8\quot SU(2),\Z)$ must be
$$(\chi(S^{-})+y^2)^2=0.$$

Since $M$ is an $S^{11}$-bundle over
$S^8\quot SU(2)$, we have one more relation in $H^{12}(M,\Z)$, saying that
the Euler class of this $S^{11}$-bundle is zero. The $S^{11}$-bundle is
$S(S^{-}\oplus V_{\R})$, and
so its Euler class is $\chi(S^{-}\oplus V_{\R})=\chi(S^{-})y$ in
$H^{12}(S^8\quot S^1,\Z)$.
Thus $M$ has cohomology ring
\begin{align*}
H^*(M,\Z)&=\Z[y,\chi(S^{-})]/((\chi(S^{-})+y^2)^2=0, \chi(S^{-})y=0)\\
&=\Z[y,\chi(S^{-})]/(\chi(S^{-})^2+y^4=0, \chi(S^{-})y=0).
\end{align*}
This is the cohomology ring of $\HH\P^4\# -\Ca\P^2$, not of $\HH\P^4\#
\Ca\P^2$.
So the biquotient $M$ is diffeomorphic to $\HH\P^4\# -\Ca\P^2$.

The proof that $\C\P^8\# -\Ca\P^2$ is diffeomorphic to a biquotient
$M$ is completely analogous: it is the quotient of the
$S^9$-bundle $S(S^{-}\oplus \R^2)$ over $S^8$ by a free $S^1$-action. Here
$S^1$ acts on $S^8$ and on the vector bundle $S^{-}$ by the homomorphism
$S^1\cong Spin(2)\subset Spin(9)$, and on $\R^2$ by the standard real
representation of $S^1$.
Thus $M$ is a biquotient of
the form $(Spin(9)\times S^9)/(Spin(8)\times S^1)$.

Finally, we come to the most surprising examples of biquotients, starting
with $\C\P^n\# \C\P^n$. The manifold
$\C\P^n$ has an orientation-reversing diffeomorphism from $n$ odd,
and so $\C\P^n\# \C\P^n$ for $n$ odd is diffeomorphic to
$\C\P^n\#-\C\P^n$. For $n$ even, we will show that
$\C\P^n\# \C\P^n$ is diffeomorphic to a biquotient $(S^3 \times
S^{2n-1})/(S^1)^2$,
for a free isometric action of $(S^1)^2$ on $S^3\times S^{2n-1}$.
The action is defined by the following homomorphism from $(S^1)^2$
to the maximal torus $(S^1)^2\times (S^1)^n$ of $SO(4)\times SO(2n)$:
$$(x,y)\mapsto ((x,y), (xy^{-1},xy,\ldots,xy)).$$
It is straightforward to check that this action of $(S^1)^2$ on $S^3\times
S^{2n-1}$
is free. Let $M$ be the quotient manifold. We use that
 the action of $(S^1)^2$ on $S^3$
has cohomogeneity one, with trivial generic stabilizer group and stabilizers
at the two
special orbits equal to the two factors $S^1$. It follows
that $M$ is the union of a 2-disc bundle
over $S^{2n-1}/S^1=\C\P^{n-1}$ (corresponding to the action of the first
factor
$S^1$) and a 2-disc bundle over $S^{2n-1}/S^1=\C\P^{n-1}$ (corresponding
to the second $S^1$) along their common boundary, which is a sphere
$S^{2n-1}$.
So $M$ is the connected sum of two copies of $\C\P^n$, with some
orientations.
To work out the orientations, it is convenient to compute the cohomology ring
of $M$.

The quotient manifold $M$ fits into a fibration
$$S^3\times S^{2n-1}\arrow M\arrow (BS^1)^2.$$
Here $(BS^1)^2$ has cohomology ring $\Z[u,v]$
with $u$ and $v$ in $H^2$.
From the description of the $(S^1)^2$-action,
we read off that the Euler classes
of the $S^3$-bundle and $S^{2n-1}$-bundle over $(BS^1)^2$ are
$uv$ and $(u-v)(u+v)^{n-1}$. Since these form a regular sequence in
the polynomial ring
$H^*((BS^1)^2,\Z)$, we have
\begin{align*}
H^*(M,\Z) &=\Z[u,v]/(uv=0, (u-v)(u+v)^{n-1}=0) \\
&= \Z[u,v]/(uv=0, u^n=v^n).
\end{align*}
This is the cohomology ring of $\C\P^n\# \C\P^n$, and (for $n$ even) not
that
of $\C\P^n\# -\C\P^n$. Therefore the biquotient $M$ is diffeomorphic
to $\C\P^n\# \C\P^n$.

We next show that $\HH\P^n\# \HH\P^n$ is diffeomorphic to a biquotient.
Imitating the construction for $\C\P^n\# \C\P^n$ might suggest identifying
$\HH\P^n\# \HH\P^n$ with a biquotient $(S^7\times S^{8n-1})/SU(2)^2$.
That works for $n$ odd, but not for $n$ even.  For $n$ even, we have to
consider
a more general type of biquotient, $(Sp(4)\times S^{4n-1})/SU(2)^3$. Here
the group $SU(2)^3$ acts on $Sp(4)$ by a homomorphism $SU(2)^3\arrow
Sp(4)^2$.
To define this,
first let $V_i$ for $i=1,2,3$ denote the standard 2-dimensional
complex representation of the $i$th factor of $SU(2)$. Then the homomorphism
$SU(2)^3\arrow Sp(4)^2$ is defined by
$(V_1\oplus V_2,V_3\oplus \C^2)$. Also, let $W_{12}:SU(2)^2\arrow SO(4)$
be the natural double covering, viewed as a 4-dimensional real
representation
of the first two copies of $SU(2)$, and let $(V_3)_{\R}:SU(2)\arrow SO(4)$
be the natural real faithful representation of the third copy of $SU(2)$. We
define the action of $SU(2)^3$ on $S^{4n-1}$ by the homomorphism
$SU(2)^3\arrow SO(8n)$ defined by $(V_3)_{\R}^{\oplus n-1}\oplus W_{12}$.

The action of $SU(2)^3$ on $Sp(4)$ has cohomogeneity one, with
trivial generic stabilizer group and with stabilizers at the two special
orbits
both isomorphic to $SU(2)$, the first conjugate to the subgroup
$\{ (x,1,x)\}$ in $SU(2)^3$ and the second conjugate to $\{(1,x,x)\}$.
These two subgroups act freely on $S^{4n-1}$, and so $SU(2)^3$ acts freely
on $Sp(4)\times S^{4n-1}$. Let $M$ be the quotient manifold. Using
the description of the $SU(2)^3$-action on $Sp(4)$, we see that $M$
is the union of two disc bundles over $S^{4n-1}/SU(2)=\HH\P^{n-1}$ along
their common boundary, $S^{4n-1}$. Therefore $M$ is the connected sum
of two copies of $\HH\P^n$, with some orientations.

To determine the relevant orientations, it suffices for $n$ even to compute
the
cohomology ring of $M$. For $n$ odd, the manifolds $\HH\P^n\#\HH\P^n$
and $\HH\P^n\# -\HH\P^n$ have isomorphic cohomology rings. In that case
we will also need a mod 3 Steenrod operation to see that $M$
is diffeomorphic to $\HH\P^n\# \HH\P^n$ rather than $\HH\P^n\# -\HH\P^n$.
First, we compute the cohomology ring of the
homotopy quotient $Sp(4)\quot SU(2)^3$, with respect to the above action
of $SU(2)^3$ on $Sp(4)$. Write $z_i=-c_2V_i$, for $i=1,2,3$, which are
generators in $H^4$ of the polynomial ring $H^*(BSU(2)^3,\Z)$.
The generators $c_2$ and $c_4$
of $H^*(BSp(4),\Z)$ pull back under the left homomorphism
$SU(2)^3\arrow Sp(4)$ to $-z_1-z_2$ and $z_1z_2$, and under the right
homomorphism $SU(2)^3\arrow Sp(4)$ to $-z_3$ and 0. Therefore,
using the Eilenberg-Moore spectral sequence as suggested
by Singhof \cite{Singhof},
the cohomology of $Sp(4)\quot SU(2)^3$ is
$$\Z[z_1,z_2,z_3]/(z_1+z_2=z_3, z_1z_2=0).$$
The manifold $M$ is an $S^{4n-1}$-bundle over $Sp(4)\quot SU(2)^3$,
up to homotopy. The Euler class of this bundle, $(V_3)_{\R}^{\oplus n-1}
\oplus W_{12}$, is $\pm (-z_3)^{n-1}(-z_1+z_2)$, and so $M$ has cohomology
ring
\begin{align*}
&=\Z[z_1,z_2,z_3]/(z_1+z_2=z_3,z_1z_2=0,z_3^{n-1}(z_1-z_2)=0) \\
&=\Z[z_1,z_2]/(z_1z_2=0, z_1^n=z_2^n).
\end{align*}
This is the cohomology ring of $\HH\P^n\# \HH\P^n$, and for $n$ even it is
not the cohomology ring of $\HH\P^n\# -\HH\P^n$. So $M$ is diffeomorphic
to $\HH\P^n\# \HH\P^n$ for $n$ even. For $n$ odd, $n\geq 3$, we also need
to observe that the classes $z_i$ are distinguished from their negatives
by the fact that $P^1z_i=-z_i^2$ in $H^8(M,\F_3)$.  By definition of the
natural
orientation on $\HH\P^n$ for $n\geq 2$, this means that $M$ is diffeomorphic
to $\HH\P^n\# \HH\P^n$ for all $n\geq 2$, not to $\HH\P^n\# -\HH\P^n$.
The case $n=1$ is trivial, since $\HH\P^n\#\HH\P^n=S^4\# S^4=S^4$.

The last Cheeger manifolds which are diffeomorphic to biquotients
are the manifolds
$\C\P^{4e+2}\# \HH\P^{2e+1}$. The relevant biquotient $M$ has the
form $M=(S^5\times S^{8e+3})/(S^1\times SU(2))$. To be more explicit,
let $L$ be the standard 1-dimensional complex representation of $S^1$,
and let $V$ be the standard 2-dimensional complex representation of $SU(2)$.
Then we let $S^1\times SU(2)$ act on $S^5$ as the unit sphere in
$V\oplus L$, and on $S^{4n-1}$ as the unit sphere
in $(V\otimes_{\C} L)^{\oplus 2e+1}$. The action of $S^1\times SU(2)$ on
$S^5$ has cohomogeneity one, with trivial generic stabilizer and with
stabilizers at the two special orbits conjugate to the two factors $S^1$ and
$SU(2)$,
respectively. Both of these subgroups act freely on $S^{8e+3}$, and so
$S^1\times SU(2)$ acts freely on $S^5\times S^{8e+3}$. Let $M$ be the
quotient manifold.

From the action of $S^1\times SU(2)$ on $S^5$, we see that $M$
is the union
of a disc bundle over $S^{8e+3}/S^1=\C\P^{4e+1}$ and a disc bundle over
$S^{8e+3}/SU(2)=\HH\P^{2e}$ along their common boundary, $S^{8e+3}$.
It follows that $M$ is the connected sum of $\C\P^{4e+2}$ and $\HH\P^{2e+1}$
with some orientations.

To show that $M$ is diffeomorphic to $\C\P^{4e+2}\#\HH\P^{2e+1}$
rather than to $\C\P^{4e+2}\# -\HH\P^{2e+1}$,
we compute the cohomology ring and
a mod 3 Steenrod operation on $M$.
Since $M=S(V\oplus L)\times S((V\otimes L)^{\oplus 2e+1})/(S^1\times SU(2))$,
we can view $M$ as an $(S^5\times S^{8e+3})$-bundle over $BS^1\times
BSU(2)$.
Here $BS^1\times BSU(2)$ has cohomology ring $\Z[x,z]$, where we let
$x=c_1L$ and $z=-c_2V$.
The vector bundles $V\oplus L$ and $(V\otimes L)^{\oplus 2e+1}$
on $BS^1\times BSU(2)$  have Euler classes
$(-z)x$ and $(x^2-z)^{2e+1}$. So $M$ has cohomology ring
\begin{align*}
&=\Z[x,z]/(-xz=0,(x^2-z)^{2e+1}=0) \\
&= \Z[x,z]/(xz=0,x^{4e+2}=z^{2e+1}).
\end{align*}
Here the element $z$ in $H^4(M,\Z)$ is distinguished
from its negative,
for $e\geq 1$, by the property that $P^1z=-z^2$ in $H^8(M,\F_3)$.
By definition of the natural orientation on $\HH\P^{2e+1}$ for $e\geq 1$,
$M$ is diffeomorphic to $\C\P^{4e+2}\# \HH\P^{2e+1}$, not to
$\C\P^{4e+2}\# -\HH\P^{2e+1}$.

By contrast, the analogous calculation of the cohomology ring
shows that the biquotient
$$(S(V\oplus L)\times S((V\otimes L)^{\oplus 2e}))/(S^1\times SU(2))$$
is diffeomorphic to $\C\P^{4e}\# -\HH\P^{2e}$, not to $\C\P^{4e}\#
\HH\P^{2e}$. This will be used in section \ref{negative}.

\section[Simplifying a given biquotient manifold]{Simplifying the description
of a given biquotient manifold}
\label{simplify}

Here is an elementary but essential
beginning to our simplification  of the description
of a given biquotient manifold. (Throughout the paper, $G$ and $H$ will
denote
compact Lie groups.)

\begin{lemma}
\label{simpconn}
Let $M$ be a simply connected biquotient manifold. Then we can write
$M=G/H$ for some simply connected group $G$ and connected group $H$
acting on $G$ by a homomorphism $H\arrow (G\times G)/Z(G)$.
If $M$ is 2-connected,
then $H$ is simply connected and $H$ acts on $G$
by a homomorphism $H\arrow G\times G$.
\end{lemma}

{\bf Proof. }Since $M$ is connected,
we can write $M$ as a biquotient $G/H$ with $G$ connected. Since $M$
is simply connected, the long exact sequence of the fibration
$H\arrow G\arrow M$,
$$\pi_1H\arrow \pi_1G\arrow \pi_1M\arrow \pi_0 H\arrow \pi_0G,$$
shows that $H$ is connected and $\pi_1H\arrow \pi_1G$ is surjective.
Let $C$ be the kernel of $\pi_1H\arrow \pi_1G$, a finitely generated abelian
group.

Let $\til{G}$ and $\til{H}$ be the universal covers of $G$ and $H$.
We can identify $\pi_1H$ with the kernel of $\til{H}\arrow H$, and
so $C$ is a central subgroup of $\til{H}$.
We have $\til{G}\cong K_G\times \R^a$ and $\til{H}\cong K_H\times \R^b$
for some simply connected compact Lie groups $K_G$ and $K_H$. The
given homomorphism $H\arrow (G\times G)/Z(G)$ lifts
to a homomorphism $\til{H}\arrow \til{G}\times \til{G}$. The resulting
action of $\til{H}$ on $\til{G}$ is trivial on the subgroup $C$, and
$\til{H}/C$ acts freely on $\til{G}$ with $M=\til{G}/(\til{H}/C)$.

Here $\til{G}=K_G\times \R^a$, and the action of $\til{H}/C$ on
$\til{G}$
is the product of an action on $K_G$ and an action on $\R^a$.
Furthermore, since the group $\R^a$ is
abelian,
$\til{H}/C$ acts on $\R^a$ by translations (left or right translations
being
the same). Since $\til{H}/C$ is a connected group and the quotient is
compact,
$\til{H}/C$ must act transitively on $\R^a$. Let $L$ be the kernel of the
action of $\til{H}/C$ on $\R^a$. Then $L$ must act freely on $K_G$,
with $M=K_G/L$. Here $K_G$ is a simply connected compact Lie group.
Also, $L$ must be connected by the long exact sequence
$$\pi_1M\arrow \pi_0L\arrow \pi_0K_G.$$

Finally, if $M$ is 2-connected, then the same long exact sequence shows that
$L$ is simply connected. In this case, the homomorphism $L\arrow (K_G\times
K_G)/Z(G)$ lifts uniquely to a homomorphism $L\arrow K_G\times K_G$.
\qed

Since we intend to study simply connected biquotient manifolds,
Lemma \ref{simpconn} tells us that we can assume $G$ is simply connected,
and thus a product of simply connected simple groups. Thus, much of the
complexity of more general compact Lie groups is avoided.

The main method of simplifying the description of a given biquotient
is the following easy observation.

\begin{lemma}
\label{remove}
Let $H$ be a compact Lie group acting on manifolds $X_1$ and $X_2$
such that $H$ acts transitively  on $X_1$ and $H$ acts freely on $X_1\times
X_2$.
Let $K\subset H$ be the stabilizer of some point in $X_1$. Then the
quotient manifold $(X_1\times X_2)/H$ is diffeomorphic to $X_2/K$,
where $K$ acts freely on $X_2$ by the restriction of the action of $H$.
\end{lemma}

{\bf Proof. }This is clear by identifying $X_1$ with $H/K$. \qed

Applying this method of simplification to biquotients gives the following
fundamental result.

\begin{lemma}
\label{trans}
Let  $M$ be a simply connected biquotient manifold. Then we can write
$M=G/H$ such that $G$ is simply connected, $H$ is connected, and $H$
does not act transitively on any simple factor of $G$.
\end{lemma}

 {\bf Proof. }By Lemma \ref{simpconn}, we can write $M$ as a biquotient
$G/H$ with
$G$ simply connected and $H$ connected. If $H$ acts transitively on some
simple
factor of $G$, then Lemma \ref{remove} allows us to remove that factor of
$G$
while replacing $H$ by a subgroup. The exact sequence $\pi_1M\arrow
\pi_0H\arrow \pi_0G$ shows that the new subgroup $H$ is still connected.
By induction on the number of simple
factors
of $G$, we can arrange that $H$ does not act transitively on any simple
factor
of $G$. \qed

\begin{convention}
\label{simp}
From now on, we will only consider
simply connected biquotient manifolds, and in writing $M=G/H$ we will
assume that $G$ and $H$ satisfy the properties listed in Lemma \ref{trans}.
\end{convention}

\section{Bounding $G$ in terms of $M$}

As mentioned in Convention \ref{simp}, for the rest of the paper we only
consider
simply connected biquotient manifolds. We can and do assume that
every such manifold $M$ is written $M=G/H$ with $G$ simply connected,
$H$ connected, and $H$ not acting transitively on any simple factor of $G$.
In this section, we will show how these properties determine $G$ up to
finitely
many possibilities in terms of the rational homotopy groups of $M$.

The following lemma will not
be needed in this section, but is included here for later use.

\begin{lemma}
\label{pi3}
Let $M=G/H$ be a simply connected biquotient manifold. If in addition
the rational homotopy group
$\pi_3M_{\Q}$ is zero, then each simple factor $H_1$ of $H$
acts trivially on all factors of $G$ isomorphic to $H_1$.
\end{lemma}

{\bf Proof. }
By the long exact sequence of the
fibration $H\arrow G\arrow M$, since $\pi_3M_{\Q}=0$,
the homomorphism $\pi_3H_{\Q}\arrow \pi_3G_{\Q}$
is surjective. In particular, for each simple factor $G_1$
of $G$,
$\pi_3H_{\Q}\arrow \pi_3(G_1)_{\Q}$
is surjective.
Suppose that some simple factor $H_1$ of $H$ which is isomorphic
to $G_1$ acts nontrivially on $G_1$;
we will derive a contradiction. Here, a priori, $H_1$ can act by left and
right
translation on $G_1$. If $H_1$ acts nontrivially on only one side of $G_1$,
thus by a nontrivial homomorphism $H_1\arrow G_1$, then this homomorphism
must be an isomorphism, since $H_1$ is simple.
In particular, $H_1$ acts transitively on $G_1$, which
contradicts Convention \ref{simp}.
Therefore, $H_1$ must act nontrivially on both sides
of $G_1$.
Then the two homomorphisms from $H_1$
to $G_1$ must both be isomorphisms, using simplicity of $H_1$ again.
We now use that any automorphism of a simply connected simple group $G_1$
acts trivially on $\pi_3(G_1)\cong \Z$, as follows from the proof of this
isomorphism
using the Killing form. The homomorphism $\pi_3(H_1)_{\Q}\arrow
\pi_3(G_1)_{\Q}$ is the difference of the homomorphisms given by the two
homomorphisms $H_1\arrow G_1$, and so it is zero. Furthermore,
since the two homomorphisms $H_1\arrow G_1$ both have only
finite centralizer, there is no room for the rest of $H$ to act on $G_1$; in
other words,
$H$ acts on $G_1$ through a quotient group isogenous to $H_1$. It follows
that the homomorphism $\pi_3H_{\Q}\arrow \pi_3(G_1)_{\Q}$ is zero,
contradicting what we know from $\pi_3M_{\Q}=0$. Thus, each simple factor
$H_1$ of $H$ must act trivially on simple factors of $G$ isomorphic
to $H_1$. \qed

To find more precise information on $G$, we need the classification
of the simple compact Lie groups.
In particular, we use that a simple group $G$ has rational
homotopy groups concentrated in odd degrees $2d-1$, and we call the numbers
$d$ that occur the {\it degrees }of $G$. Particularly important for us is
the
{\it maximal degree }of $G$, which we call $d(G)$. It is also called
the {\it Coxeter number }of $G$.
The degrees of $G$ are well known in many
contexts: we can also say that $H^*(G,\Q)$ is an exterior algebra with
generators
in degrees $2d-1$ where $d$ runs over the degrees of $G$, or that the
degrees of $G$
are the degrees of the generators of the ring of invariants of the Weyl
group acting
on its reflection representation. We tabulate the degrees of the simple
groups here, following Bourbaki \cite{Bourbaki} or Gorbatsevich-Onishchik
(\cite{GO}, Table 1, p.~127).
To avoid repetitions, one can assume that $A_l$ has $l\geq 1$, $B_l$ has
$l\geq 3$,
$C_l$ has $l\geq 2$, and $D_l$ has $l\geq 4$.

\begin{tab}
\label{degrees}
\begin{align*}
A_l &:  2,3,\ldots,l+1\\
B_l &: 2,4,6,\ldots,2l\\
C_l &: 2,4,6,\ldots,2l\\
D_l &: 2,4,6,\ldots,2l-2;l\\
G_2 &: 2,6\\
F_4 &: 2,6,8,12\\
E_6 &: 2,5,6,8,9,12\\
E_7 &: 2,6,8,10,12,14,18\\
E_8 &: 2,8,12,14,18,20,24,30
\end{align*}
\end{tab}

Using Table \ref{degrees}
and the known low-dimensional representations of each
group,
Onishchik proved the following result \cite{Onishchikinclusion}. He later
gave a
more systematic proof, using reflection groups \cite{Onishchikreflection}.

\begin{lemma}
\label{max}
Let $H\arrow G$ be a nontrivial
homomorphism of simply connected simple groups. Then the maximal degrees
satisfy $d(H)\leq d(G)$. Moreover, if $d(H)=d(G)$, then either $H\arrow G$
is an isomorphism or $G/H$ is one of the following homogeneous spaces.
On the right we show
the degrees of $G$ and $H$.
$$\begin{array}{lrl}
Spin(2n)/Spin(2n-1)=S^{2n-1},
  n\geq 4& 2,4,6,\ldots,2n-2;n &2,4,6,\ldots,2n-2 \\
SU(2n)/Sp(2n), n\geq 2& 2,3,4,\ldots,2n & 2,4,6,\ldots 2n \\
Spin(7)/G_2=S^7 &2,4,6 & 2,6 \\
Spin(8)/G_2=S^7\times S^7 & 2,4,4,6 & 2,6\\
E_6/F_4& 2,5,6,8,9,12& 2,6,8,12
\end{array}$$

In all these cases except $Spin(8)/Spin(7)=S^7$, there is a unique
conjugacy class
of nontrivial homomorphisms $H\arrow G$; in the case $Spin(8)/Spin(7)$,
there are three conjugacy classes which are equivalent under outer
automorphisms
of $Spin(8)$. Also, in all the above cases, the centralizer of $H$ in $G$ is
finite.
\end{lemma}

The following result shows how to apply Lemma \ref{max} to biquotients
$M=G/H$,
although it is only a step on the way to the more precise Theorem
\ref{class}.
For a simple factor
$G_1$ of $G$, we say that a degree $d$ of $G_1$ is {\it killed by $H$ }if
the homomorphism $\pi_{2d-1}H_{\Q}\arrow \pi_{2d-1}(G_1)_{\Q}$
associated to the action of $H$ on $G_1$ is nonzero.

\begin{lemma}
\label{maxapp}
Let $M=G/H$ be a simply connected biquotient,
written according to Convention \ref{simp}.
 Let $G_1$ be a simple factor
of $G$ such that the maximal degree of $G_1$ is killed by $H$. Then either
there is a simple factor $H_1$ of $H$ such that $H_1$ acts
nontrivially on exactly one side of $G_1$ by one of the homomorphisms
in Lemma \ref{max}, so that $G_1/H_1$ is one of
$Spin(2n)/Spin(2n-1)=S^{2n-1}$,
$SU(2n)/Sp(2n)$, $Spin(7)/G_2=S^7$, $Spin(8)/G_2=S^7\times S^7$,
or $E_6/F_4$;
or $G_1$ is isomorphic to $SU(2n+1)$ for some $n$ and there is a simple
factor $H_1$ of $H$ also isomorphic to $SU(2n+1)$ which acts on $G_1$
by $h(g)=hgh^{t}$. The $SU(2n+1)$ case cannot occur if $\pi_3M_{\Q}=0$.
\end{lemma}

{\bf Proof. }We are given a simple factor $G_1$ of $G$ such that the maximal
degree of $G_1$ is killed by $H$. It follows that there must be a simple
factor
$H_1$ of $H$ which kills the maximal degree $d$ of $G_1$.
Since $H_1$ is simply connected,
the action of $H_1$ on $G_1$ is given by a homomorphism
$H_1\arrow G_1\times G_1$, and the resulting linear map
$\pi_{2d-1}(H_1)_{\Q}\arrow \pi_{2d-1}(G_1)_{\Q}\cong \Q$
is nonzero. This linear map is the difference of the two linear maps
associated to the two homomorphisms $H_1\arrow G_1$ (on the left and right),
so at least one of those two linear maps is nonzero. By Lemma \ref{max},
either $H_1$ is isomorphic to $G_1$ or $(G_1,H_1)$ is one of the pairs
$(Spin(2n),Spin(2n-1))$, $(SU(2n),Sp(2n))$, $(Spin(7),G_2)$,
$(Spin(8),G_2)$,
or $(E_6,F_4)$.

Suppose first that $(G_1,H_1)$ is one of these 5 pairs.
If $H_1$ acts nontrivially
on
both sides of $G_1$, then in all cases except $(Spin(8),Spin(7))$,
Lemma \ref{max} implies that the two
homomorphisms $H_1\arrow G_1$ are conjugate, so the resulting map
$\pi_*(H_1)_{\Q}\arrow \pi_*(G_1)_{\Q}$ (the difference of the left and
right
maps)
is 0, contradicting  the fact that $H_1$ kills the top degree of $G_1$. Even
in the
case $(Spin(8),Spin(7))$, we compute that
the outer automorphism group of $Spin(8)$ acts trivially on the top degree,
6,
of $Spin(8)$, that is, on $\pi_{11}Spin(8)_{\Q}$. It follows
that if $Spin(7)$
acts nontrivially on both sides of $Spin(8)$, then the left and right
homomorphisms
$Spin(7)\arrow Spin(8)$ give the same linear map into the top degree of
$Spin(8)$,
and so the action of $Spin(7)$ of $Spin(8)$ cannot kill the top degree of
$Spin(8)$,
a contradiction. Thus, in all these cases, $H_1$ acts nontrivially on only
one side
of $G_1$, by one of the homomorphisms in Lemma \ref{max}. The lemma is
proved in this case.

The remaining case is where $H_1$ is isomorphic to $G_1$. Clearly this
cannot
occur if $\pi_3M_{\Q}=0$,  by Lemma \ref{pi3}. In general,
if $H_1$ acts
nontrivially on only one side of $G_1$, then it must act by an isomorphism
$H_1\arrow G_1$. So $H$ acts transitively on $G_1$, contrary to
Convention \ref{simp}.
Therefore $H_1$ must act nontrivially on both sides of $G_1$, clearly by
two isomorphisms $H_1\arrow G_1$. We are given that $H_1$ kills the top
degree
of $G_1$, so these two isomorphisms $H_1\arrow G_1$ must give
{\it different }linear maps into the top degree of $\pi_*(G_1)_{\Q}$. So
the outer automorphism group of $G_1$ must act nontrivially on the top
degree
of $G_1$. Of the simply connected simple groups, only $A_n$, $D_n$, and
$E_6$ have nontrivial outer automorphism group, and only in the case
$G_1=SU(2n+1)=A_{2n}$ does the outer automorphism group act nontrivially
on the top degree of $G_1$. (The outer automorphism group $\Z/2$ of $SU(n)$,
$n\geq 3$, acts by the identity on the even degrees of $SU(n)$ and by $-1$ on
the
odd degrees. A topological way to see this is to identify $\pi_*G_{\Q}$
with
the dual to the vector space $H^{>0}(BG,\Q)/(H^{>0}\cdot H^{>0})$,
and use in the case $G=SU(n)$ that the outer automorphism $E\mapsto E^*$
acts on Chern classes by $c_i\mapsto (-1)^ic_i$.)
So we must have $G_1=SU(2n+1)$, with
$H_1=SU(2n+1)$ acting on $G_1$ by the identity on one side and by
the outer automorphism $x\mapsto (x^t)^{-1}$ on the other. \qed

In order to have strong restrictions on the simple factors of $G$ in terms
of the rational homotopy groups of $M=G/H$, we need to analyze more
completely the way
$H$ acts on simple factors of $G$ isomorphic to $Spin(2n)$, $SU(2a)$,
$Spin(7)$, $Spin(8)$, or $E_6$. In analogy with Lemmas \ref{max} and
\ref{maxapp}, we will first formulate a general statement on the
subgroups
of these groups, and then apply it to biquotients. Before that, we need some
simple topological results.

\begin{lemma}
\label{surjlemma}
Let $H$ and $G$ be connected compact
Lie groups, and let $f:H\arrow G$ be any
continuous map. If $f_*:\pi_*H_{\Q}\arrow \pi_*G_{\Q}$ is surjective,
then $f$ is surjective.
\end{lemma}

{\bf Proof. }Both $H$ and $G$ have the rational homotopy type of products
of odd-dimensional spheres. So $H^*(H,\Q)$ and $H^*(G,\Q)$ are both
exterior algebras, and the assumption means that the homomorphism
$$f^*:H^{>0}(G,\Q)/(H^{>0}\cdot H^{>0})\arrow H^{>0}(H,\Q)/
(H^{>0}\cdot H^{>0})$$
is injective. Thus $H^*(H,\Q)$ is the exterior algebra generated by the
generators of $H^*(G,\Q)$ together with some other generators. So
$H^*(G,\Q)\arrow H^*(H,\Q)$ is injective. Since $G$ is a closed orientable
manifold,
it follows that $H\arrow G$ is surjective. \qed

\begin{corollary}
\label{surj}
Let $H$ and $G$ be connected compact Lie groups, with an action of $H$
on $G$ by a homomorphism $H\arrow (G\times G)/Z(G)$. If the associated
homomorphism $\pi_*H_{\Q}\arrow \pi_*G_{\Q}$ is surjective, then
$H$ acts transitively on $G$.
\end{corollary}

{\bf Proof. }We are given that the homomorphism associated to the orbit map
$f:H\arrow G$ of some point in $G$ is surjective on rational homotopy
groups.
By Lemma \ref{surjlemma}, $H\arrow G$ is surjective. That is, $H$ acts
transitively on $G$. \qed

We now apply Corollary \ref{surj} to get information on the subgroups
of the groups occurring in Lemma \ref{max}. (We also use the classification
of simple Lie groups in the following proof, but it is more pleasant to use
Corollary \ref{surj} when possible.)

\begin{lemma}
\label{secondmax}
Let $\varphi:H\arrow Spin(2n)$ be a homomorphism from a simply connected
simple
group such that the homomorphism
$$\pi_{2n-1}H_{\Q}\arrow \pi_{2n-1}Spin(2n)_{\Q}/\pi_{2n-1}Spin(2n-1)_{\Q}
\cong \Q$$
is not zero. Then $H$ acts transitively on the sphere
$Spin(2n)/Spin(2n-1)=S^{2n-1}$.

Next, let $G/K$ be one of the other homogeneous spaces from Lemma \ref{max}:
$SU(2n)/Sp(2n)$ with $n\geq 2$, $Spin(7)/G_2=S^7$,
$Spin(8)/G_2=S^7\times S^7$, or $E_6/F_4$. Let $H\arrow G$
be a homomorphism from a simply connected compact Lie group $H$
which kills the second-largest degree of $G$. (In the
case $G/K=Spin(8)/G_2$, where $Spin(8)$ has degrees $2,4,4,6$, we assume
that both degrees 4 of $Spin(8)$ are killed by $H$.)
Then $H$ acts transitively on $G/K$.
\end{lemma}

{\bf Proof. }In the cases where $G/K$ is one of
$Spin(2n)/Spin(2n-1)=S^{2n-1}$, $Spin(7)/G_2=S^7$,
or $Spin(8)/G_2=S^7\times S^7$, the assumption implies that
the action of $H\times K$ on $G$ is surjective on rational homotopy groups.
By Corollary \ref{surj}, $H\times K$ acts transitively on $G$. Equivalently,
$H$ acts transitively on $G/K$.

Next, let $G/K=SU(2n)/Sp(2n)$, $n\geq 2$,
and suppose that $H$ kills the second-largest
degree, $2n-1$, of $SU(2n)$. We can replace $H$ by one of its simple factors
without changing this property. If $H$ is isomorphic to $G$,
then the homomorphism $H\arrow G$ is an isomorphism, and so $H$ acts
transitively on $G/K$. Otherwise,
by Lemma \ref{max}, $H$ has maximal degree at most that of $SU(2n)$,
which is $2n$, and if equality holds then $H=Sp(2n)$.
But $Sp(2n)$
does not kill the degree $2n-1$ of $SU(2n)$. So $H$ must have maximal degree
$2n-1$. By Table \ref{degrees},
the only simple group with maximal degree an odd number is $H=SU(2n-1)$.
Since $n\geq 2$, any nontrivial homomorphism
$SU(2n-1)\arrow SU(2n)$ is equivalent to the standard inclusion by some
automorphism of $SU(2n-1)$. Then $SU(2n-1)$ acts transitively on
$SU(2n)/Sp(2n)$,
because $Sp(2n)$ acts transitively on $SU(2n)/SU(2n-1)=S^{4n-1}$.

Finally, suppose that $G/K=E_6/F_4$ and $H$ kills the second-largest degree,
9, of
$E_6$. We can replace $H$ by one of its simple factors without changing
this property.
Then $H$ has rank at most the rank 6 of $E_6$ and has 9 as a degree.
By Table \ref{degrees}, it follows that $H$ is isomorphic
to
$E_6$. So the given homomorphism $H\arrow E_6$ must be an isomorphism,
and in particular $H$ acts transitively on $E_6/F_4$. \qed

We now apply this result on subgroups to deduce strong information on
the classification of biquotient manifolds. For a given biquotient
$M=G/H$, we say that a given simple factor $G_i$ of $G$
{\it contributes degree $d_i$ to $M$ }if the homomorphism
$\pi_{2d_i-1}H_{\Q}\arrow \pi_{2d_i-1}(G_i)_{\Q}$ is not surjective.
If $G=\prod G_i$ and every factor $G_i$ contributes some degree $d_i$
to $M$, then $\pi_{2d-1}M_{\Q}=\pi_{2d-1}G_{\Q}/\pi_{2d-1}H_{\Q}$
has dimension at least equal to the number of simple factors $G_i$ with
$d_i=d$.

\begin{theorem}
\label{class}
Let $M=G/H$ be a simply connected biquotient manifold, written using
Convention \ref{simp}. Let $G_1$ be any simple factor of $G$.
Then at least one of the following holds.

(1) $G_1$ contributes its maximal degree to $M$.

(2) $G_1$ contributes its second-largest degree to $M$, and there is a
simple
factor $H_1$ of $G_1$ which acts nontrivially on exactly one side of $G_1$,
with $G_1/H_1$ isomorphic to one of the homogeneous spaces
$SU(2n)/Sp(2n)$ with $n\geq 2$,
$Spin(7)/G_2=S^7$,
$Spin(8)/G_2=S^7\times S^7$, or $E_6/F_4$. The second-largest degree of
$G_1$ is, respectively, $2a-1$, $4$, $4$, or $9$. (In the case
$G_1=Spin(8)$,
which has degrees $2,4,4,6$, the claim is only that $G_1$ contributes at
least
one degree 4 to $M$.)

(3) $G_1\cong Spin(2n)$ with $n\geq 4$
contributes its degree $n$ to $M$, and there is a
simple
factor $H_1\cong Spin(2n-1)$ which acts nontrivially on exactly one side
of $G_1$, by the standard inclusion, with $Spin(2n)/Spin(2n-1)=S^{2n-1}$.

(4) $G_1\cong SU(2n+1)$ contributes degrees $2,4,6,\ldots,2n$
 to $M$, and there is a simple
factor $H_1\cong SU(2n+1)$ of $G_1$ which acts on $G_1$ by
$h(g)=hgh^{t}$.
\end{theorem}

{\bf Proof. }Let $G_1$ be a simple factor of $G$. Suppose that (1) does not
hold,
in other words that the maximal degree of $G_1$ is killed by $H$. By Lemma
\ref{maxapp}, there is a simple factor $H_1$ of $H$ such that either
$H_1$ acts nontrivially on exactly one side of $G_1$ by one of the
homomorphisms
listed in (2) or (3) above, or $G_1=H_1=SU(2n+1)$ and $H_1$ acts on $G_1$
as in (4). The remaining point is to show that $G_1$ contributes the degrees
to $M$ that we have claimed.

In cases (2) and (3), $H_1$ has finite centralizer in $G_1$ by Lemma
\ref{max},
so the rest of $H$ can act on $G_1$ only on the other side from $H_1$.
Since $H_1$ and the rest of $H$ together do not act transitively on $G_1$,
by Convention \ref{simp}, Lemma \ref{secondmax} shows that $G_1$ contributes
its degree $n$ to $M$ if $G_1/H_1=Spin(2n)/Spin(2n-1)=S^{2n-1}$, or its
second-largest degree to $M$ if $G_1/H_1$ is one of $SU(2n)/Sp(2n)$,
$Spin(7)/G_2=S^7$, $Spin(8)/G_2=S^7\times S^7$, or $E_6/F_4$.

In case (4), since $H_1=SU(2n+1)$ acts with finite centralizer on both sides
 of $G_1=SU(2n+1)$,
no other factor of $H$ can act on $G_1$. So the image of $\pi_*H_{\Q}
\arrow \pi_*(G_1)_{\Q}$ is equal to the image of $\pi_*(H_1)_{\Q}\arrow
\pi_*(G_1)_{\Q}$. This homomorphism is the difference of the identity map
on $\pi_*(G_1)_{\Q}$ with the map given by the outer automorphism
$g\mapsto (g^t)^{-1}$, which acts by 1 on the even degrees and by $-1$
on the odd degrees (as shown in the proof of Lemma \ref{maxapp}).
So the image of $\pi_{2d-1}(H_1)_{\Q}\arrow \pi_{2d-1}(G_1)_{\Q}$
is zero for $d$ even. That is, $G_1$ contributes all its even degrees
$2,4,\ldots,2n$ to $M$. \qed

Theorem \ref{class} implies the following important qualitative statement
on the classification of biquotients. The analogous result for homogeneous
spaces is easy and probably well known.
It is perhaps surprising
that the following statement requires the detailed classification work we
have done
in this section, but that seems to be true, at least for now.

\begin{theorem}
\label{finite}
There are only finitely many diffeomorphism classes of 2-connected
biquotient
manifolds of a given dimension.
\end{theorem}

Theorem \ref{finite}
is vaguely reminiscent of the Petrunin-Tuschmann theorem,
which says in particular that for any number $C$,
there are only finitely many diffeomorphism classes
of 2-connected closed Riemannian manifolds with curvature
$0\leq K\leq C$ and diameter 1 \cite{PT}. But there is probably no
way to actually deduce Theorem \ref{finite} from
the Petrunin-Tuschmann theorem, since there is no obvious upper
bound on the curvature of 2-connected biquotients until one goes through the
proof of Theorem \ref{finite}. In fact, the discussion after
Theorem 1.1 in my paper \cite{Totarocurv}
shows that there can be no upper bound on the curvature of simply
connected biquotients of dimension 6, if one fixes their diameter to be 1.

{\bf Proof of Theorem \ref{finite}. }Any
biquotient manifold $M$ is rationally elliptic;
that is, all but finitely many of the rational homotopy groups of $M$
are zero. So,
writing $n$ for the dimension of $M$, the odd-degree rational homotopy
groups $\pi_{2d-1}M_{\Q}$  are zero for $d>n$, and the total dimension
of the odd-degree rational homotopy groups of $M$ is at most $n$,
by Friedlander and Halperin (\cite{FHTbook}, p.~434).
Since $M$ is simply connected,
we can write $M$ as a biquotient $G/H$
according to Convention \ref{simp}.
In particular, $G$ is simply connected and,
since $M$ is 2-connected, $H$ is also simply connected. The essential
point is Theorem
\ref{class}, which implies that
each simple factor of $G$ contributes at least one degree to $M$,
and that one such degree
is at least half the maximal degree of $G$. Therefore
the number of simple factors of $G$ is at most
the total dimension of $\pi_{\odd}M_{\Q}$ and hence at most $n$,
and the maximal degree of each simple factor of $G$
is at most $2n$. Thus $G$, being a product of simply connected simple
groups, is determined up to finitely many possibilities by $n$. Also, $H$ is
a simply connected group of dimension at most that of $G$, so $H$ and
the homomorphism $H\arrow (G\times G)/Z(G)$ (up to conjugacy) are
determined up to finitely many possibilities. \qed

\section[Further general results on classification]{Further general results
on the classification of biquotients}

In this section we continue the previous section's method: we classify
subgroups
of compact Lie groups with certain properties, and apply the results to the
classification of biquotients. The main result of this section is
Theorem \ref{top}, which describes the possible simple factors of $G$
in a biquotient $M=G/H$ which contribute only their top degree to $M$.
(This terminology is defined before Theorem \ref{class}.)

We begin by stating a classification of certain subgroups of $Spin(2n)$.

\begin{lemma}
\label{spin}
Let $\varphi:H\arrow Spin(2n)$, $n\geq 4$, be a homomorphism from a simply
connected simple group $H$ such that the linear map
$$\pi_{2n-1}H_{\Q}\arrow \pi_{2n-1}Spin(2n)_{\Q}/\pi_{2n-1}Spin(2n-1)_{\Q}
\cong \Q$$
is not zero. Then $\varphi$ is either an isomorphism or one of the following
homomorphisms, up to the standard $\Z/2$ group of outer automorphisms
of $Spin(2n)$: $SU(n)\inj Spin(2n)$, $Sp(2a)\inj Spin(4a)$, the spin
representation
$Spin(7)\inj Spin(8)$, or the spin representation $Spin(9)\inj Spin(16)$.
\end{lemma}

This is straightforward to prove using the known degrees and
low-dimensional representations of all the simple groups. Alternatively,
we can deduce it from Borel's
classification
of groups which act linearly and transitively on the sphere, in the
odd-dimensional
case \cite{Borel}. Namely, the hypothesis of Lemma \ref{spin} means that
the action of $H$ on $S^{2n-1}$ associated to the homomorphism
$H\arrow Spin(2n)$ has orbit map $H\arrow S^{2n-1}$ such that
$\pi_{2n-1}H_{\Q}\arrow \pi_{2n-1}S^{2n-1}_{\Q}\cong \Q$ is not zero.
It follows that $H\arrow S^{2n-1}$ is surjective, in other words that
$H$ acts transitively on $S^{2n-1}$. Then Lemma \ref{spin} follows
from Borel's classification.

The subgroups we classify next are the simple subgroups $H$ of any
simple group $G$ such
that the maximal degree of $H$ is at least the second-largest degree of $G$.
This is closely related to two classifications by Onishchik. First, he
classified
the simple subgroups
$H\subset G$ such that
$\dim_{\Q}\pi_{\odd}(G/H)_{\Q}=1$ (\cite{GO}, Table 3, p.~185),
that is, $H$ kills all but one degree of $G$; these make up the first part
of the list  in Lemma \ref{second}. Next, he classified the simple subgroups
$H\subset G$
such that $d(H)\geq d(G)-2$ (\cite{Onishchikbook}, Table 7.2, p.~195); it
turns out
that these include all the subgroups on the second
part of the list in Lemma \ref{second}.

\begin{lemma}
\label{second}
Let $H\arrow G$ be a nontrivial homomorphism of simply connected simple
groups such that the maximal degree of $H$ is at least
the second-largest degree of
$G$ and is less than the maximal degree of $G$. Then 
$G/H$ is isomorphic to one of the following homogeneous spaces.
On the left is the Dynkin index of the homomorphism $H\arrow G$. On the
right are shown the degrees of $G$ not occurring in $H$, and the degrees of
$H$
not occurring in $G$, in both cases with multiplicities. In the last column
is the centralizer of $H$ in $G$, written modulo finite groups.
$$\begin{array}{rlrlc}
1 & SU(n)/SU(n-1)=S^{2n-1}, n\geq 3 & n &  & S^1\\
1 & Sp(2n)/Sp(2n-2)=S^{4n-1}, n\geq 2 & 2n &  & A_1\\
1 &Spin(2n+1)/Spin(2n)=S^{2n},n\geq 3 & 2n & n & 1\\
1 & Spin(2n+1)/Spin(2n-1)=UT(S^{2n}), n\geq 3 & 2n & & S^1\\
2 & Sp(4)/SU(2)=UT(S^4) & 4 & & S^1\\
10 & Sp(4)/SU(2) & 4 & & 1\\
4 & SU(3)/SO(3) & 3 & & 1\\
1 & Spin(9)/Spin(7)=S^{15} & 8 & & 1\\
1 & G_2/SU(3)=S^6 & 6 & 3 & 1\\
1 & G_2/SU(2)=UT(S^6) & 6 & & A_1\\
3 & G_2/SU(2) & 6 & & A_1\\
4 & G_2/SO(3) & 6 & & 1\\
28 & G_2/SO(3) & 6 & & 1\\
1 & F_4/Spin(9)=\Ca\P^2 & 12 & 4 & 1\\
\\
1 & Spin(2n)/Spin(2n-2)=UT(S^{2n-1}), n\geq 4 & n,2n-2 & n-1 & S^1\\
1 & Spin(2n)/Spin(2n-3), n\geq 4& n,2n-2 & & A_1\\
1 & SU(2n+1)/Sp(2n), n\geq 2 & 3,5,\ldots,2n+1 & & S^1\\
2 & SU(2n+1)/SO(2n+1), n\geq 2& 3,5,\ldots,2n+1 & & 1\\
1 & Spin(10)/Spin(7) & 5,8 & & S^1\\
2 & SU(7)/G_2 & 3,4,5,7 & & 1\\
1 & Spin(9)/G_2 & 4,8 &  &S^1\\
1 & Spin(10)/G_2 & 4,5,8 &  &A_1
\end{array}$$
\end{lemma}

The proof is straightforward for $G$ classical, using the known
low-dimensional
representations
of each simple group, or alternatively the results by Onishchik mentioned
above.
For $G$ exceptional, more than enough information is provided by
Dynkin's paper
on the subgroups of the exceptional groups \cite{Dynkin}. For example,
Dynkin's Table 16 classifies the $A_1$ subgroups of $G_2$.
Notice that the listings for $Spin(9)/Spin(7)=S^{15}$ and $Spin(10)/Spin(7)$
refer
to the spin representation of $Spin(7)$; these spaces are different from
the spaces $Spin(2n+1)/Spin(2n-1)=UT(S^{2n})$ for $n=4$ and
$Spin(2n)/Spin(2n-3)$ for $n=5$. Also, there are three conjugacy classes of
nontrivial
homomorphisms $SU(2)\arrow Sp(4)$, where we write $V$ for the standard
representation of $SU(2)$: $V\oplus \C^2$, where
$Sp(4)/SU(2)=Sp(4)/Sp(2)=S^7$, which has Dynkin index 1;
 $V\oplus V$, where $Sp(4)/SU(2)=Spin(5)/Spin(3)=UT(S^4)$, which has
Dynkin index 2; and $S^3V$, which has Dynkin index 10.

We now apply Lemma \ref{second} to the classification of biquotients.
Together with Theorem \ref{class}, the following theorem will be
our most important tool in the classification of biquotients.

\begin{theorem}
\label{top}
Let $M=G/H$ be a simply connected biquotient, written according
to Convention \ref{simp}. Let $G_1$ be a simple factor
of $G$ which contributes only its top degree to $M$. Then at least one
of the following holds.

(1) $G_1$ is isomorphic to $SU(2)$, which has degree 2.

(2) $G_1$ is a rank-2 group $SU(3)$,
$Sp(4)$, or $G_2$, with top degree 3, 4, 6 respectively, and there is a
simple
factor $H_1\cong SU(2)$ of $H$ which acts nontrivially on $G_1$.

(3) There
is a simple factor $H_1$ of $H$ such that $H_1$ acts nontrivially on exactly
one side of $G_1$ and $G_1/H_1$ is one of the homogeneous spaces
in the first part of Lemma \ref{second}'s list.

(4) There are two simple factors $H_1$ and $H_2$
of $H$ which act nontrivially
on the two sides of $G_1$ in one of the following ways,
up to switching $H_1$ and $H_2$:
$G_1=Spin(2n)$, $n\geq 4$, $H_1$ is $Spin(2n-2)$ or $Spin(2n-3)$,
$H_2$ is $SU(n)$, or $Sp(2a)$ with $n=2a$, or $Spin(9)$ with $n=8$;
here $G_1$ has top degree $2n-2$. Or $G_1=SU(2n+1)$, $n\geq 3$,
$H_1$ is $Sp(2n)$ or $SO(2n+1)$, and $H_2=SU(2n-1)$; here $G_1$
has top degree $2n+1$. Or $H_1\back G_1/H_2$ is one of
$G_2\back SU(7)/SU(5)$, $G_2\back Spin(9)/SU(4)$,
$G_2\back Spin(9)/Sp(4)$, or $G_2\back Spin(10)/SU(5)$; here
$G_1$ has top degree 7, 8, 8, 8, respectively.
\end{theorem}

{\bf Proof. }If $G_1=SU(2)$, we have conclusion (1). So we can assume
that $G_1$ has rank at least 2. Equivalently, $G_1$ has
at least two degrees. There must
be a simple factor $H_1$ of $H$ which kills at least one
second-largest degree of $G_1$. (By Table \ref{degrees},
$G_1$ has a unique second-largest degree except when
$G_1$ is $Spin(8)$, which has degrees $2,4,4,6$.)

Suppose that $H_1$ is isomorphic to $G_1$. If $H_1$ acts nontrivially on
only one side of $G_1$, then it acts by an isomorphism on that side of
$G_1$,
and so $H_1$ acts transitively on $G_1$, contrary to Convention \ref{simp}.
So $H_1$ must act by an isomorphism on both sides of $G_1$. Since
all automorphisms of $G_1$ act as the identity on $\pi_3G_1=\Z$, the
resulting
homomorphism $\pi_3(H_1)_{\Q}\arrow \pi_3(G_1)_{\Q}$ is zero.
Also, since $H_1$ is acting with finite centralizer on both sides of $G_1$,
the rest of $H$ cannot act on $G_1$, and so the whole homomorphism
$\pi_3H_{\Q}\arrow \pi_3(G_1)_{\Q}$ is zero. That is, $G_1$ contributes
its degree 2 to $M$. Since $G_1$ is not isomorphic to $SU(2)$,
this contradicts our
assumption that $G_1$ contributes only its top degree to $M$.

Thus $H_1$ is not isomorphic to $G_1$.
We know that $d(H_1)\leq d(G_1)$ by Lemma \ref{max}. Suppose
that $d(H_1)=d(G_1)$. By Lemma \ref{max}, $(G_1,H_1)$ is one of the pairs
$(Spin(2n),Spin(2n-1)$,
$(SU(2n),Sp(2n))$, $(Spin(7),G_2)$, $(Spin(8),G_2)$, or $(E_6,F_4)$.
In all these cases except $(Spin(8),Spin(7))$,
there is a unique conjugacy class of nontrivial homomorphisms
$H_1\arrow G_1$; moreover, even in the case $(Spin(8),Spin(7))$,
all nontrivial homomorphisms are equivalent under automorphisms of
$Spin(8)$ and so they all give the same homomorphism $\pi_3H_1\arrow
\pi_3G_1$.
Also, the centralizer of $H_1$ in $G_1$ is finite in all cases.
So, in all cases, if $H_1$ acts nontrivially on both sides of $G_1$, then
no other factor of $H$ acts on $G_1$, and the homomorphism
$\pi_3H_{\Q}\arrow \pi_3(G_1)_{\Q}$ is zero. That is,
$G_1$ contributes its degree 2 to $M$, contrary to our assumption
that $G_1$ contributes only its top degree to $M$. Therefore $H_1$ must
act only on one side of $G_1$. But that implies, in all these cases, that
$H_1$ kills the top degree of $G_1$. This contradicts our assumption
that $G_1$ contributes its top degree to $M$.

So we must have $d(H_1)<d(G_1)$. Since $H_1$ kills at least one
second-largest degree of $G_1$,
$d(H_1)$ must be at least the second-largest degree
of $G_1$. Therefore $(G_1,H_1)$ must be one of the pairs listed in
Lemma \ref{second}. If $H_1$ is isomorphic to $SU(2)$, then $G_1$
has rank 2 and we have conclusion (2).
We can now assume that $H_1$ is not isomorphic to $SU(2)$.

Next, we will show that $H_1$ acts nontrivially on only one side of $G_1$
in all the remaining cases. Suppose that $H_1$ acts nontrivially on both
sides of $G_1$. By Lemma \ref{second}, the centralizer of $H_1$ on each
side of $G_1$ is at most finite by $A_1$. So any simple factor of $H$
other than $H_1$ which acts nontrivially on $G_1$ is isomorphic
to $SU(2)$. Thus $H_1$ by itself must kill all the degrees of $G_1$
greater than 2 and less than the maximal degree $d(G_1)$. This is clearly
impossible for the pairs $(G_1,H_1)$ on the second part of Lemma
\ref{second}'s list, since in these cases $G_1$ contains at least one
degree in the interval $(2,d(G_1))$ with greater multiplicity than
$H_1$ does.

So $(G_1,H_1)$ is on the first part of Lemma \ref{second}'s list. Since
$H_1$ is not isomorphic to $SU(2)$, the list shows that all nontrivial
homomorphisms $H_1\arrow G_1$ have Dynkin index 1. Since $H_1$ acts
nontrivially on both sides of $G_1$, it  follows that the associated
homomorphism
$\pi_3(H_1)_{\Q}\arrow \pi_3(G_1)_{\Q}$ is zero. That is, $H_1$ does
not kill the degree 2 of $G_1$. Since the whole group $H$ does kill
the degree 2 of $G_1$, at least one of the two nontrivial
homomorphisms $H_1\arrow G_1$ must have centralizer containing
an $A_1$ subgroup. By the first part of Lemma \ref{second}'s list,
it follows that the pair $(G_1,H_1)$ is $(Sp(2n),Sp(2n-2))$ for some $n\geq 3$.

But for $(G_1,H_1)$ equal to $(Sp(2n),Sp(2n-2))$
with $n\geq 3$, there is a unique
conjugacy class of nontrivial homomorphisms $H_1\arrow G_1$. Since
$H_1$ acts nontrivially on both sides of $G_1$, it follows that
the associated homomorphism $\pi_*(H_1)_{\Q}\arrow \pi_*(G_1)_{\Q}$
is zero. In particular, $H_1$ does not kill the second-largest
degree, $2n-2$, of $G_1$. This is a contradiction.
We have thus completed the proof
that $H_1$ acts nontrivially on only one side of
$G_1$. 

For spaces $G_1/H_1$ on the first part of Lemma \ref{second}'s list,
this completes the proof of conclusion (3).
It remains to consider spaces $G_1/H_1$ on the second part of Lemma
\ref{second}'s list. To prove conclusion (4), we need to identify
another simple factor of $H$ which acts on $G_1$.

In several cases on the second part of Lemma \ref{second}'s list,
$G_1$ is isomorphic to $Spin(2n)$, $n\geq 4$. Here $H_1$ is either
$Spin(2n-2)$,
$Spin(2n-3)$, or, for $n=5$, $Spin(7)$ (with a different homomorphism to
$Spin(10)$)
or $G_2$.
In all these cases, the degree
$n$ of $G_1$ is not killed by $H_1$, meaning that
$\pi_{2n-1}(H_1)_{\Q}\arrow \pi_{2n-1}(G_1)_{\Q}$ is not surjective.
So there must be another simple factor $H_2$ of $H$ which acts
on $G_1$ such that $\pi_{2n-1}(H_2)_{\Q}\arrow \pi_{2n-1}(G_1)_{\Q}/
\pi_{2n-1}(H_1)_{\Q}$ is not zero. If $n$ is odd, this just means
that $\pi_{2n-1}(H_2)_{\Q}\arrow \pi_{2n-1}(G_1)_{\Q}$ is not zero.
On the other hand, if $n$ is even, then we can always assume, after
automorphisms of $H_1$ and $G_1$, that $H_1\arrow G_1$ is the standard
inclusion. So, for $n$ odd or even, we can say that
$$\pi_{2n-1}(H_2)_{\Q}\arrow \pi_{2n-1}Spin(2n)_{\Q}/\pi_{2n-1}Spin(2n-1)
_{\Q}\cong \Q$$
is not zero.
Here $H_2$ cannot be $SU(2)$ since
$n\geq 4$,
so $H_2$ must act only on the other side of $G_1$ from $H_1$, since
the centralizer of $H_1$ in $G_1$ is at most finite by $A_1$. By Lemma
\ref{spin},
the space $G_1/H_2$ must be one of  $Spin(2n)/SU(n)$, $Spin(4a)/Sp(2a)$,
$Spin(8)/Spin(7)$, or $Spin(16)/Spin(9)$.  If $G_1/H_2=Spin(8)/Spin(7)$,
then we can replace $H_1$ by $H_2$ and we have conclusion (3); so we can
exclude that case. We have thus proved conclusion (4) for $G_1=Spin(2n)$.

Next, there are several cases in the second part of Lemma \ref{second}'s
list where $G_1=SU(2n+1)$, $n\geq 2$. Here $H_1$ is either $Sp(2n)$,
$SO(2n+1)$, or, for $n=3$, $G_2$.
In all these cases, $H_1$ has only even
degrees, and in particular it does not kill the degree $2n-1$ of $G_1$. So
there
must be another simple factor $H_2$ of $H$ which kills the degree $2n-1$
of $G_1$. Here $H_2$ must act nontrivially only on
the other side of $G_1$ from $H_1$, because the centralizer of $H_1$ in
$G_1$
is at most $S^1$ by finite in these cases. It is easy to read from Table
\ref{degrees} that any simple group $H_2$ which maps
nontrivially to $G_1=SU(2n+1)$, $n\geq 2$,  and has $2n-1$ as a degree is
either
$SU(2n-1)$, $SU(2n)$, or $SU(2n+1)$. Here $H_2=SU(2n+1)$ is excluded
by Convention \ref{simp}, which says
that $H$ does not act transitively on $G_1$. If
$H_2=SU(2n)$, then we can replace $H_1$ by $H_2$ and we have conclusion
(3). The remaining possibility is $H_2=SU(2n-1)$. We have proved conclusion
(4) for $G_1=SU(2n+1)$.

The last case, from the second part of Lemma \ref{second}'s list, is
where $G_1$ is $Spin(9)$ and $H_1$ is the exceptional group $G_2$.
Here the degree 4 of $G_1$ is not killed by $H_1$, so it must be
killed by some other simple factor $H_2$ of $H$. The centralizer of $H_1$
in $G_1$ is finite by $S^1$, so $H_2$ must act nontrivially only on the
other
side of $G_1$ from $H_1$. From Table \ref{degrees},
since $H_2$ has a degree 4 and maps nontrivially to $Spin(9)$, $H_2$ is one
of $Sp(4)$, $SU(4)$, $Spin(7)$, $Spin(8)$, or $Spin(9)$. The case
$H_2=Spin(9)$
is excluded by Convention \ref{simp}, which says
that $H$ does not act transitively on $G_1$.
If $H_2$ is $Spin(7)$ or $Spin(8)$, then we can replace $H_1$ by $H_2$ and
we have conclusion (3). The remaining possibilities for $H_2$
are $Sp(4)$ and $SU(4)$.
We have proved conclusion (4) for $G_1=Spin(9)$. \qed

\section[Rational homology spheres]{Biquotients which are rational homology
spheres}

As an application of Theorems \ref{class} and \ref{top}, we now classify
all biquotients which are simply connected rational homology spheres.
The result seems surprising. In particular, the Gromoll-Meyer exotic
7-sphere which is a biquotient \cite{GM} is the only such example in any
dimension. As mentioned in the introduction, Theorem \ref{sphere} was
found at the same time by Kapovitch and Ziller \cite{KZ}.

\begin{theorem}
\label{sphere}
Any biquotient which is a simply connected rational homology sphere
is either a homogeneous manifold, the Gromoll-Meyer exotic 7-sphere
which is a biquotient $Sp(4)/SU(2)$, or a certain
4-connected 11-manifold with the integral homology groups of $UT(S^6)$
which is a biquotient $G_2/SU(2)$.

The homogeneous manifolds which are simply connected rational homology
spheres are the sphere $S^n$, the unit tangent
bundle
$UT(S^{2n})$, the Wu 5-manifold $SU(3)/SO(3)$ \cite{Dold},
the Berger 7-manifold
$Sp(4)/SU(2)$ with $\pi_3M$ isomorphic to $\Z/10$ \cite{Berger},
the 11-manifold $G_2/SU(2)$
with $\pi_3M$ isomorphic to $\Z/3$, and the two 11-manifolds $G_2/SO(3)$
with
$\pi_3M$ isomorphic to $\Z/4$ or $\Z/28$.
\end{theorem}

The homogeneous spaces in Theorem \ref{sphere} were classified
by Onishchik (\cite{GO}, Table 3, p.~185).
The nontrivial biquotient $G_2/SU(2)$ was first discovered by
Eschenburg (\cite{EschenburgHab}, pp.~166--170).
Kapovitch and Ziller
have shown that this biquotient $G_2/SU(2)$ is diffeomorphic
to the connected sum of $UT(S^6)$ with some homotopy 11-sphere
\cite{KZ}. It is not known whether $G_2/SU(2)$ is actually
diffeomorphic to $UT(S^6)$; one hopes for a negative answer,
which would be more interesting.

Before starting the proof of Theorem \ref{sphere},
we assemble some elementary facts about biquotients $Sp(4)/H$. In the
following
lemma, we consider free actions of a group $H$ on $Sp(4)$ given by a
homomorphism
$H\arrow Sp(4)^2/Z(Sp(4))$.

\begin{lemma}
\label{sp4}
(1) There is no free action of $SO(3)$ on $Sp(4)$.
Any free action of $SU(2)$ on $Sp(4)$ is either trivial on one side of
$Sp(4)$,
so that $Sp(4)/SU(2)$ is $S^7$, $UT(S^4)$, or the Berger 7-manifold,
or given by the homomorphisms $(V\oplus \C^2,V^{\oplus
2})$
up to switching the two sides of $Sp(4)$,
so that $Sp(4)/SU(2)$ is the Gromoll-Meyer exotic 7-sphere \cite{GM}.
Here $V$ denotes
the standard representation of $SU(2)$.

(2) Any free action of $SU(2)^2$ on $Sp(4)$ is given, up to switching the
two
$SU(2)$ factors and switching the two sides of $Sp(4)$, by the homomorphisms
$(V_1\oplus V_2,\C^4)$ or $(V_1\oplus\C^2,V_2^{\oplus 2})$. Here
$V_1$ and $V_2$ are the standard representations of the two $SU(2)$ factors.
In both cases, the quotient $Sp(4)/SU(2)^2$ is diffeomorphic to $S^4$.
\end{lemma}

{\bf Proof. }We only prove (1) here, the calculation for (2) being
similar.
There are three conjugacy classes of nontrivial
homomorphisms $SU(2)\arrow Sp(4)$, $V\oplus\C^2$, $V^{\oplus 2}$,
and $S^3V$, by Lemma \ref{second}.
If $SU(2)$ acts trivially on one side of $Sp(4)$, then
$Sp(4)/SU(2)$
is one of the three homogeneous spaces mentioned in (1).
So we can assume that $SU(2)$ acts nontrivially on both sides
of $Sp(4)$. Let $H$ be $SU(2)$ or $SU(2)/\{ \pm 1\}\cong SO(3)$.
The group $H$ acts
freely on $Sp(4)$ if and only if the two images of each nontrivial element
of $H$
in $Sp(4)$ are not conjugate. In particular,
the two homomorphisms $SU(2)\arrow Sp(4)$
must be non-conjugate. If one homomorphism is $V\oplus \C^2$ and the other
is $V^{\oplus 2}$, then Gromoll and Meyer showed that $SU(2)$ acts freely
on $Sp(4)$ and that the quotient manifold is an exotic 7-sphere \cite{GM}.
Otherwise, $SU(2)$ must act by $S^3V$ on one side and either $V\oplus \C^2$
or $V^{\oplus 2}$ on the other.
Then neither $SU(2)$ nor $SU(2)/\{\pm 1\}$
acts freely on $Sp(4)$,
since the images of the diagonal matrix $(\zeta_3,\zeta_3^{-1})$ in $SU(2)$
under $S^3V$ and $V\oplus \C^2$ are conjugate, and likewise the images
of the diagonal matrix $(\zeta_4,\zeta_4^{-1})$ in $SU(2)$ under $S^3V$
and $V^{\oplus 2}$ are conjugate. (Here $\zeta_n$ denotes a primitive $n$th
root of unity.) This proves (1). \qed

{\bf Proof of Theorem \ref{sphere}. }We
 write $M=G/H$ according to Convention \ref{simp}.
Since $M$ is a rational homology sphere, the odd-dimensional rational
homotopy of $M$ has dimension 1. In more detail, $S^{2n-1}$ has
1-dimensional
rational homotopy in dimension $2n-1$ and zero otherwise, while $S^{2n}$
has 1-dimensional rational homotopy in dimensions $2n$ and $4n-1$, and zero
otherwise.
Since each simple factor of $G$ contributes
at least dimension 1 to $\pi_{\odd}M_{\Q}$ by Theorem \ref{class}, $G$
must be simple. This already makes the situation much more understandable;
it would not be clear at all without the analysis leading to Theorem
\ref{class}.

If the simply connected rational homology sphere $M$ has dimension $r$ at
most 4,
then it is automatically a homotopy sphere. It is therefore not surprising
to find that
$M$ is diffeomorphic to $S^r$ for $r\leq 4$. First, for $r=2$ we can just
use that
every homotopy 2-sphere is diffeomorphic to $S^2$. For $r\geq 3$, we have
$\pi_2M_{\Q}=0$, and so $H$ has finite fundamental group;
equivalently, $H$ is semisimple. For $r=3$,
$\pi_3M\cong \Z$ and $\pi_4M_{\Q}=0$, which implies that $G$ has one more
simple factor than $H$ does. Since $G$ is simple, $H=1$.
Since $M$ is a homotopy 3-sphere, $G$ has degree 2
only,
and so $G$ is $SU(2)$. Thus $M$ is diffeomorphic to $SU(2)$, that is,
to $S^3$.

For $r=4$, $\pi_2M=\pi_3M=0$
 and $\pi_4M\cong \Z$, which implies that $H$ is simply connected and has one
more
simple factor than $G$ does. Since $G$ is simple, $H$ has two simple
factors.
By the rational homotopy groups of $M$, $G$ must contribute degree 4 to $M$
while
$H$ contributes degree 2, and nothing else. By Theorems \ref{class} and
\ref{top},
there is a simple factor $H_1$ of $H$ such that either $G/H_1$ is a
homogeneous
space diffeomorphic to $S^7$, or $G/H_1$ is $Spin(8)/G_2=S^7\times S^7$,
or $(G,H_1)$ is $(Sp(4),SU(2))$.
Let $H_2$ denote the other simple factor of $H$.
If $G/H_1=S^7$, then $H_2$ has degree 2 only,
so $H_2$ is isomorphic to $SU(2)$. In this case $M$ is the quotient of $S^7$
by
$SU(2)$ acting freely by a homomorphism $SU(2)\arrow O(8)$. Such a
homomorphism
is unique up to conjugacy, and so $M$ is the standard quotient $S^4$. Next,
if $G/H_1$ is $Spin(8)/G_2=S^7\times S^7$,
then $H_2$ has degrees 2 and 4, and so $H_2$
is isomorphic to $Sp(4)$. But there is in fact no free isometric action of
$Sp(4)$
on $S^7\times S^7$,  because
the restriction of any action of $Sp(4)$ on $S^7$ to the subgroup
$SU(2)=Sp(2)\subset Sp(4)$ has a fixed point. (Either the associated
8-dimensional complex representation of $Sp(4)$ has a trivial summand,
or it is the sum of two copies of the standard representation of
$Sp(4)$ and hence restricts on $SU(2)$ to $V^{\oplus 2}\oplus \C^4$.)

For $r=4$, it remains to consider the case $(G,H_1)=(Sp(4),SU(2))$. Here
$H_2$
has degree 2 only, and so $H_2$ is isomorphic to $SU(2)$. Thus $M$ is a
biquotient
$Sp(4)/SU(2)^2$. By Lemma \ref{sp4},
$M$ is diffeomorphic to $S^4$.

We proceed to the case $r\geq 5$. As mentioned earlier, we have
$\pi_2M_{\Q}=0$, and so $H$ is semisimple.
Since $\pi_3M_{\Q}=\pi_4M_{\Q}=0$, $H$ has the same number of simple
factors as $G$. Since $G$ is simple, so is $H$.

One of the cases (1) to (4) in Theorem \ref{class} must hold. Here (4) is
excluded since $G$ is simple and $H$ acts freely on $G$.
In case (3), $M$ is the homogeneous
space $Spin(2n)/Spin(2n-1)=S^{2n-1}$. In case (2), $M$ is one of the
homogeneous
spaces $SU(2a)/Sp(2a)$ with $a\geq 2$, $Spin(7)/G_2=S^7$, $Spin(8)/G_2
=S^7\times S^7$, or $E_6/F_4$. Since $M$ is a rational homology sphere,
considering the degrees of these homogeneous spaces shows that $M$ is
either $SU(4)/Sp(4)=Spin(6)/Spin(5)=S^5$ or $Spin(7)/G_2=S^7$.

There remains case (1) of Theorem \ref{class}, where $G$ contributes its
top degree to $M$. In this case, one of conclusions (1) to (4) in Theorem
\ref{top}
must hold. 
Here case (4) of Theorem \ref{top} is excluded since $H$ is simple. In case
(3),
$M$ is one of the homogeneous spaces listed in the first part of Lemma
\ref{second}. Since $M$ is a rational homology sphere, the possibilities
are as listed in Theorem \ref{sphere}.

Case (1), $G=SU(2)$, of Theorem \ref{top} is excluded because we are
considering
biquotients $M$ of dimension $r\geq 5$. So there remains only case (2). That
is, $H_1$ is $SU(2)$,  $H$ is either $H_1=SU(2)$
or $H_1/\{ \pm 1\}\cong SO(3)$, and
$M$ is a biquotient $SU(3)/H$, $Sp(4)/H$, or $G_2/H$, of dimension
5, 7, or 11, respectively.
The homogeneous spaces of this type are listed in Lemma \ref{second}.
So we can assume that $H_1$ acts nontrivially on both sides of $G$.
Biquotients of this type (a rank-2 group divided by a rank-1 group)
were classified by Eschenburg (\cite{EschenburgHab}, pp.~166--170),
but I will give
my own proof since Eschenburg's paper is not widely available.

First, suppose $G=SU(3)$. Then there are two conjugacy classes of nontrivial
homomorphisms $SU(2)\arrow SU(3)$, $V\oplus \C$ and $S^2V$, where $V$
denotes the standard representation of $SU(2)$. Since $H$ (which is $SU(2)$
or $SO(3)$) acts freely on $G$,
the two images in $SU(3)$ of any nontrivial element of $H$ are not
conjugate in $SU(3)$.
So $SU(2)$
must act by $V\oplus \C$ on one side of $SU(3)$ and by $S^2V$ on the other.
But, if we
write $\zeta_a$ for a primitive $a$th root of unity, the diagonal matrix
$(\zeta_3, \zeta_3^{-1})$ in $SU(2)$ has images under both homomorphisms to
$SU(3)$
which are conjugate to $(\zeta_3,1,\zeta_3^{-1})$. So neither $SU(2)$
nor $SU(2)/\{ \pm 1\}$ can be
acting freely. Thus, the case $G=SU(3)$ does not occur.

The case $G=Sp(4)$ is covered by Lemma \ref{sp4}. Since $SU(2)$ acts
nontrivially on both sides of $Sp(4)$, $M$ is the Gromoll-Meyer exotic
7-sphere.

The last case to consider is where $G$ is the exceptional group
$G_2$. Here, as described in Table
\ref{second},
there are 4 conjugacy classes of nontrivial homomorphisms $SU(2)\arrow G_2$,
which we identify by their Dynkin index, 1, 3, 4, or 28. Using these Dynkin
indices,
we compute that the composed representation $SU(2)\arrow G_2\arrow U(7)$
must be: $W_1=V^{\oplus 2}\oplus \C^3$, $W_3=S^2V\oplus V^{\oplus 2}$,
$W_4=(S^2V)^{\oplus 2}\oplus \C$, and $W_{28}=S^6V$. In particular,
the homomorphisms $W_4$ and $W_{28}$ from $SU(2)$ to $G_2$
are trivial on $\{ \pm 1\}$ in
$SU(2)$, and the other two are not.

We are assuming that $SU(2)$ acts nontrivially on both sides of $G_2$. Since
the action is free for either $SU(2)$ or $SU(2)/\{ \pm 1\}$, the two
homomorphisms
$SU(2)\arrow G_2$ must be non-conjugate. There are now 6 cases to consider,
corresponding to the 6 unordered pairs of distinct homomorphisms $W_i$.
Notice that the action of $SU(2)$ on $G_2$ is trivial on $\{ \pm 1\}$
if and only if the two homomorphisms send $-1$ to the same
element
of the center of $G_2$. The center of $G_2$
is trivial, so the two homomorphisms
must both
be trivial on $-1$. That is, $SU(2)$ acts on $G_2$ through its
quotient $SO(3)$ if and only if the two homomorphisms are $(W_4,W_{28})$.

Given that $SU(2)$ acts on $G_2$ by $W_i$ on one side and $W_j$ on the
other,
the action of $SU(2)$
is free if and only if every nontrivial element of $SU(2)$ has
non-conjugate
images in $G_2$ under the two homomorphisms. Every element of $SU(2)$
is conjugate to a diagonal element $(x,x^{-1})$ with $x\in S^1$, so it
suffices
to consider those elements. Furthermore, it is convenient to observe that
two elements
of $G_2$ are conjugate if and only if their images in $U(7)$ are conjugate.

Using this, we can check whether each pair $(W_i,W_j)$ of homomorphisms
$SU(2)\arrow G_2$ gives a free action of $SU(2)$ on $G_2$ or not. 
Here the image of $x\in S^1\subset SU(2)$ in $U(7)$ under the
homomorphism $W_i$ is conjugate to the diagonal matrix:
\begin{align*}
W_1 &: (x,x^{-1},x,x^{-1},1,1,1) \\
W_3 &: (x,x^{-1},x,x^{-1},x^2,1,x^{-2}) \\
W_4 &: (x^2,x^{-2},x^2,x^{-2},1,1,1) \\
W_{28} &: (x^{-6},x^{-4},x^{-2},1,x^2,x^4,x^6).
\end{align*}
The result is that 5 of the 6 unordered pairs $(W_i,W_j)$ do not
give free actions of $SU(2)$ (or $SO(3)$, in the case $(W_4,W_{28})$)
on $G_2$. Indeed, the following nontrivial elements $x\in S^1$
have conjugate images (in $U(7)$, hence in $G_2$) under the
two representations.
\begin{align*}
(W_1,W_3) &: x=-1\\
(W_1,W_4) &: x=\zeta_3\\
(W_1,W_{28}) &: x=\zeta_3\\
(W_3,W_{28}) &: x=\zeta_5\\
(W_4,W_{28}) &: x=\zeta_3
\end{align*}

If $SU(2)$ acts on $G_2$ by $(W_3,W_4)$, however, then we compute
from the above formulas that the two images of any nontrivial element
$x\in S^1\subset SU(2)$ are not conjugate in $U(7)$, and hence
not conjugate in $G_2$. So this is a free
action of $SU(2)$ on $G_2$. Since the Dynkin indices 3 and 4 differ by 1,
the resulting biquotient $M=G_2/SU(2)$ has $\pi_3M=0$. Using the
known homology of $G_2$, we compute that the 11-manifold
$M$ is 4-connected and has the integral homology groups of $UT(S^6)$.

Thus, the only biquotient $G_2/H$ with $H$ isomorphic to
$SU(2)$ or $SO(3)$ and $H$ acting
nontrivially on both sides is the 4-connected 11-manifold $G_2/SU(2)$,
with $SU(2)$ acting by $(W_3,W_4)$. \qed

\section[Three Cheeger manifolds]{Three Cheeger manifolds which are not
homotopy equivalent to biquotients}
\label{16}

We now return to the proof of Theorem \ref{cheeger}. Roughly
in order of increasing difficulty, we show in this section
that the 16-manifolds
$\Ca\P^2\# \Ca\P^2$, $\C\P^8\# \Ca\P^2$, and $\HH\P^4\# \Ca\P^2$
are not homotopy equivalent to biquotients.

Suppose that the 16-manifold
$\Ca\P^2 \# \Ca\P^2$ is homotopy equivalent to
a biquotient $M=G/H$. As throughout the paper, we assume that
$M$ is written as a biquotient $G/H$ which satisfies Convention \ref{simp}.
Here $H^*(M,\Z)\cong \Z[x,y]/(xy=0,x^2=y^2)$,
$|x|=|y|=8$. This is a complete intersection ring with 2 generators
in degree 8 and 2 relations in
degree 16. It follows that the rational homotopy groups of $M$
are isomorphic to $\Q^2$ in dimensions 8 and 15, and otherwise 0 (the same
as for $S^8\times S^8$).
So the map $\pi_{2d-1}H_{\Q}\arrow \pi_{2d-1}G_{\Q}$ must have
2-dimensional cokernel for $d=8$, 2-dimensional kernel for $d=4$,
and otherwise is an isomorphism. Equivalently, we say that $G$ contributes
two degrees 8 to $M$, while $H$ contributes two degrees 4, and nothing else.

Each simple factor of $G$ contributes at least one degree to $M$, by
Theorem \ref{class}, so $G$ has at most two simple factors. Suppose
first that $G$ is simple. We know that $G$ contributes two degrees 8 to $M$,
and nothing else. One of the four cases in Theorem \ref{class}
must hold. Here (1) cannot hold: if $G$ contributes its maximal degree to
$M$,
that would have to be 8, but the maximal degree of each simple Lie group
occurs only with multiplicity 1. The remaining cases of the theorem are
incompatible
with the fact that $G$ contributes only degree 8 to $M$.
Thus we have a contradiction from the assumption
that $G$ is simple.

So $G$ has two simple factors, $G=G_1\times G_2$. (Here $G_2$
does not denote the exceptional group $G_2$.) Each factor must contribute
one degree 8 to $M$, and nothing else. We can apply Theorem \ref{class}
to each factor $G_i$. From the degrees, it is clear that only cases (1)
and (3) can occur. That is, for $1\leq i\leq 2$, either $G_i$ contributes
its maximal degree to $M$, which must be 8, or else $G_i=Spin(16)$
and there is a simple factor $H_i\cong Spin(15)$ of $H$ which acts
nontrivially on exactly one side of $G_i$, with $G_i/H_i=S^{15}$.
If $G_i$ contributes its maximal degree to $M$, we can apply Theorem
\ref{top} to $G_i$. Cases (1) and (2) cannot arise, since the maximal degree
of $G_i$ is 8. Therefore, by cases (3) and (4),
 there is a simple factor $H_i$ of $H$ which
acts nontrivially on exactly one side of $G_i$ such that $G_i/H_i$
is one of the homogeneous spaces $Spin(9)/Spin(8)=S^8$, 
$SU(8)/SU(7)=Sp(8)/Sp(6)=Spin(9)
/Spin(7)=S^{15}$,
$Spin(9)/Spin(7)=UT(S^8)$,
$Spin(10)/Spin(8)$, $Spin(10)/Spin(7)$, $Spin(9)/G_2$, or
$Spin(10)/G_2$.

Since neither $G$ nor $H$ contributes any degree 2 to $M$, and each
simple group has exactly one degree 2 (that is, every simple group
has $\pi_3=\Z$), $H$ must have the same number
of simple factors as $G$. Thus $H$ has two simple factors.

Let $G_1$ and $G_2$ denote the two simple factors of $G$. We know
that there is a simple factor $H_1$ of $H$ such that $H_1$ acts
nontrivially on exactly one side of $G_1$, with $G_1/H_1$ equal to
either $Spin(16)/Spin(15)=S^{15}$ or one of the other homogeneous
spaces listed above. We also know the analogous statement for $G_2$,
but a priori it is possible that the same simple factor $H_1$ of $H$
plays the same role for both $G_1$ and $G_2$.
 But then $G=G_1\times G_2$ has two
degrees 14 (if $G_1$ and $G_2$ are $Spin(16)$) or 7 (if $G_1$
and $G_2$ are $SU(8)$) or
6 (if $G_1$ and $G_2$ are $Spin(9)$ or $Spin(10)$),
 both of which must be killed by $H$, whereas the group $H_1$ has only
one degree 14 or 6 or 7. So we can order the two simple
factors of $H$ in such a way that $H_i$ kills the relevant degree of $G_i$,
both for $i=1$ and for $i=2$.
By the proof of Theorem \ref{top}, it follows that $H_i$ acts nontrivially
on exactly one side of $G_i$, with $G_i/H_i$ equal to one of the
above homogeneous spaces, for $i=1$ and for $i=2$.

Since $M$ has dimension only 16, both $G_1/H_1$ and $G_2/H_2$
must be the homogeneous space $Spin(9)/Spin(8)=S^8$. If $H_1$ acts trivially
on $G_2$, or also if $H_2$ acts trivially on $G_1$, then $M$ is an
$S^8$-bundle
over $S^8$ and hence has signature zero, contradicting the fact that
$M$ is homotopy equivalent to $\Ca\P^2\# \Ca\P^2$. So $H_1$ acts
nontrivially on $G_2$ and $H_2$ acts nontrivially on $G_1$.
Since $H_i$ has finite centralizer in $G_i$ for $i=1,2$, the action of $H_1$
on $G_2$ must be given by a nontrivial homomorphism $H_1\arrow G_2$
on the other side of $G_2$ from $H_2$, and likewise for the action of
$H_2$ on $G_1$. Any nontrivial homomorphism $Spin(8)\arrow Spin(9)$
has Dynkin index 1. So the homomorphism $\pi_3H\arrow \pi_3G$,
$\Z^2\arrow \Z^2$, is given by a $2\times 2$ matrix with both rows
equal to $(1,-1)$ or $(-1,1)$.
Such a matrix has zero determinant, and so $\pi_3M_{\Q}$
is not zero. Again, this contradicts the fact that $M$ is homotopy
equivalent
to $\Ca\P^2\# \Ca\P^2$. This completes the proof that $\Ca\P^2\# \Ca\P^2$
is not homotopy equivalent to a biquotient. In fact, the proof shows
that $\Ca\P^2\# \Ca\P^2$ is not even rationally homotopy equivalent
to a biquotient.

We now prove that $\C\P^8\# \Ca\P^2$ is not homotopy equivalent to
a biquotient $M=G/K$. Since $G$ is simply connected, the boundary map
$\Z\cong \pi_2M\arrow \pi_1K$ in the long exact sequence is an isomorphism.

So, if we let $H$ be the commutator subgroup of $K$, then $H$ is simply
connected and $K/H$ is isomorphic to $S^1$. Let $N$ be the 17-manifold
$G/H$. Since $N$ is the natural $S^1$-bundle over $M\simeq \C\P^8
\# \Ca\P^2$, we compute that $N$ has the integral cohomology ring
of $S^8\times S^9$. To go further, we use Wu's theorem that
the Stiefel-Whitney classes of the tangent bundle of $M$
are invariants of the homotopy type of $M$ (\cite{MS}, Theorem 11.14).
Because
the 8-sphere in the Cayley plane has self-intersection $1$,
the Stiefel-Whitney class $w_8(\Ca\P^2)\in H^8(\Ca\P^2,\F_2)$
 is not zero. Therefore $w_8(\C\P^8\# \Ca\P^2)$
is not in the subgroup $(H^2)^4$ of $H^8$. By Wu's theorem, the same
holds for the manifold $M$ which is homotopy equivalent to
$\C\P^8\# \Ca\P^2$. Since $N$ is an $S^1$-bundle over $M$, it follows
that $w_8(N)$ is not zero. We will use this later.

The cohomology ring of $M$ is a complete intersection ring.
The degrees of the generators and relations determine the
rational homotopy groups of $M$, and hence of $N$.
The result is that $G$ contributes degrees 5 and 8 to $N$, while
$H$ contributes 4, and nothing else. Since every simple factor of $G$
contributes at least one degree to $N$, $G$ has at most 2 simple factors.
Also, since neither $G$ nor $H$ contributes degree 2, $G$ and $H$
have the same number of simple factors.

Suppose first that $G$ is simple. It follows that $H$ is simple, too.
We can apply Theorem \ref{class}, and cases (2), (3), (4) are excluded
since $G$ contributes degrees 5 and 8 to $N$ and no more. Therefore case
(1) must hold; that is, $G$ contributes its maximal degree to $N$.
So $G$ has maximal degree 8. There is no simple group whose highest
two degrees are 5 and 8, so the second-largest degree of $G$ must
be killed by $H$. Since $G$ has exactly one degree 8, $H$ must
have maximal degree less than 8.
So, regardless of how $H$ acts on $G$, $(G,H)$ must be one of the pairs in
Lemma \ref{second}. From the list there, since we know that $G$
adds degrees 5 and 8 to $M$ and $H$ adds degree 4 and no more,
we must have $G=Spin(10)$ and $H=Spin(8)$. All nontrivial homomorphisms
$H\arrow G$ have Dynkin index 1. So if $H$ acts nontrivially on both sides
of $G$, then $G$ would contribute degree 2 to $M$, a contradiction.
So $H$ acts on only one side of $G$ and $G/H$ is isomorphic to
the homogeneous space $Spin(10)/Spin(8)=UT(S^9)$.

Indeed, there is a free $S^1$-action on $UT(S^9)$, with quotient the
8-dimensional complex quadric $Q^8_{\C}$. This quadric has the rational
homotopy
type of $\C\P^8\# \Ca\P^2$. But no $S^1$-quotient of $N=UT(S^9)$ can have
the homotopy type of $\C\P^8\# \Ca\P^2$.
The point is that $N$ is an $S^8$-bundle
over $S^9$. It follows that there is a sphere $S^8$ in $N$ which represents
a generator of $H_8(N,\F_2)$ and which has trivial normal bundle.
Therefore $w_8(N)|S^8=0$, and hence $w_8(N)=0$, contradicting our
earlier calculation that the Stiefel-Whitney class $w_8(N)$ is not zero.

So $G$ must have two simple factors. Since $G$ and $H$ have
the same number of simple factors, $H$ also has two simple
factors. Each simple factor of $G$ must contribute
exactly one degree to $N$; we can assume that
$G_1$ contributes degree 8 and $G_2$
contributes degree 5. By Theorems \ref{class} and \ref{top},
there is a simple factor $H_1$ acting on exactly one side of $G_1$
such that $G_1/H_1$ is one of the homogeneous spaces
$Spin(9)/Spin(8)=S^8$, 
$Spin(16)/Spin(15)=SU(8)/SU(7)=Sp(8)/Sp(6)=Spin(9)/Spin(7)=S^{15}$,
$Spin(9)/Spin(7)=UT(S^8)$, 
$Spin(10)/Spin(8)$, $Spin(10)/Spin(7)$, $Spin(9)/G_2$, or $Spin(10)/G_2$.
Also, by the same theorems,
there is a simple factor $H_2$ acting on exactly one side
of $G_2$ such that $G_2/H_2$ is one of the homogeneous spaces
$Spin(10)/Spin(9)=SU(5)/SU(4)=S^9$ or $SU(6)/Sp(6)$.

If $H_1$ and $H_2$ are the same factor of $H$, this factor must be
isomorphic to $Sp(6)$, and $N$ is a biquotient of the form
$(Sp(8)\times SU(6))/(Sp(6)\times X)$. Here $X$ must be a simple group
with degrees $2,3,4,4,6$, but there is no such group. So $H_1$ and
$H_2$ are the two different simple factors of $H$.

From the degrees of $N$ (or just by its dimension, which is
only 17), the homogeneous space $G_1/H_1$ must be $Spin(9)/Spin(8)=S^8$
and the homogeneous space $G_2/H_2$ must be either
$Spin(10)/Spin(9)=S^9$ or $SU(5)/SU(4)=S^9$. The 17-manifold $N$
is not an $S^8$-bundle over $S^9$, because we know that $w_8(N)$
is not zero. So $H_1=Spin(8)$ must act nontrivially on the
second factor $G_2$ of $G$. It follows that $G_2/H_2$ is
$Spin(10)/Spin(9)=S^9$, not $SU(5)/SU(4)=S^9$. Furthermore,
$H_2=Spin(9)$ must act trivially on $G_1=Spin(9)$, by Lemma
\ref{pi3}.
Thus $N$ is a biquotient of the form $(Spin(9)\times Spin(10))/
(Spin(8)\times Spin(9))$ with $Spin(9)$ acting trivially on $Spin(9)$
and $Spin(8)$ acting nontrivially on $Spin(10)$. That is, $N$ is a
nontrivial
$S^9$-bundle over $S^8$.

There are three conjugacy classes of nontrivial homomorphisms from
$Spin(8)$ to $Spin(10)$. They have the form $W\oplus \R^2$,
where $W$ is either
the standard 8-dimensional
real representation $V$
of $Spin(8)$ or one of the two spin representations $S^{-}$ and $S^{+}$.
Thus $N$ is the $S^9$-bundle
$S(W\oplus \R^2)$ over $S^8$, where $W$ is the real vector bundle
over $Spin(9)/Spin(8)=S^8$ corresponding to the representation $W$
of $Spin(8)$. If $W$ is the standard representation $V$ of $Spin(8)$,
then the corresponding vector bundle on $S^8$ is the tangent bundle,
which has zero Stiefel-Whitney classes. It follows that the sphere
bundle $S(V\oplus \R^2)$ has zero Stiefel-Whitney classes, contrary
to what we know about $N$. So $W$ must be one of the two spin
representations
of $Spin(8)$. Without loss of generality, we can assume that $W=S^{-}$;
if instead $W=S^{+}$, we can apply an automorphism of order 2 of $Spin(8)$
which switches the isomorphism classes of the representations $S^{-}$
and $S^{+}$ and does not change the isomorphism class of the representation
$V$ and hence of the standard inclusion $Spin(8)\arrow Spin(9)$. Once
that is done, $N$ is the $S^9$-bundle $S(S^{-}\oplus \R^2)$ over $S^8$.

Finally, we need to analyze the free $S^1$-action on $N$ that gives
$M=N/S^1$. By computing the centralizers of the homomorphisms that
define $N$, we see that the $S^1$-action on $N$ is defined
by a homomorphism $S^1\arrow Spin(9)$ together with a homomorphism
from $S^1$ to $S^1$, the identity component of the centralizer
of $Spin(8)$ in $Spin(10)$. The homomorphism $S^1\arrow Spin(9)$
defines an action of $S^1$ on $S^8$ and also on the
$Spin(9)$-equivariant vector bundle $S^{-}$ over $S^9$, while
the homomorphism $S^1\arrow S^1$, of the form $z\mapsto z^b$ for
some integer $b$, defines the action of $S^1$ on the trivial bundle
$\R^2$ over $S^8$. Since $S^1$ is acting freely on $N=S(S^{-}\oplus \R^2)$,
it must act freely on both $S(S^{-})=S^{15}$ and on $S(\R^2)=S^1\times S^8$.
Since the $S^1$-action on $S^8$ must have a fixed point, the freeness
of the $S^1$-action on $S(\R^2)$ means that $b=\pm 1$. We compute
that there is a unique conjugacy class of homomorphisms $S^1\arrow Spin(9)$
which give a free action of $S^1$ on $S(S^{-})=S^{15}$,
namely the subgroup $S^1=Spin(2)\subset Spin(9)$.
We can assume that $b=1$, since changing $b=-1$ to $b=1$ clearly does
not change the diffeomorphism class of the quotient manifold
$M=N/S^1$. Thus we have uniquely described the biquotient $M$.

But the biquotient $M$ we have described is exactly the one which
we proved to be diffeomorphic
to $\C\P^8\# -\Ca\P^2$ in section \ref{positive}. In particular, it
is not homotopy equivalent to $\C\P^8\# \Ca\P^2$. This completes
the proof that $\C\P^8\#\Ca\P^2$ is not homotopy equivalent to a
biquotient.

We now show that $\HH\P^4\#\Ca \P^2$ is not homotopy equivalent to
a biquotient $M=G/H$. Since $M$ is 2-connected, $H$ is simply connected,
by Lemma \ref{simpconn}. The ring $H^*(M,\Z)$ is a complete intersection,
and so the degrees of its
generators and relations determine the rational homotopy
groups of $M$. The result is that $G$ contributes degrees 6 and 8 to $M$,
while $H$ contributes degrees 2 and 4, and no others. Since every simple
factor of $G$ contributes at least one degree to $M$, $G$ has at most
2 simple factors. Since $H$ contributes exactly one degree 2, $H$ has
one more simple factor than $G$ has.

Suppose first that $G$ is simple. We can apply Theorem \ref{class},
and cases (2), (3), (4) are excluded because $G$ contributes degrees
6 and 8 to $M$ and no more. Therefore case (1) must hold; that is,
$G$ contributes its maximal degree to $M$. So $G$ has maximal degree 8.
By Table \ref{degrees}, $G$ must be one of
$Sp(8)$, $Spin(9)$, $Spin(10)$, or $SU(8)$, which have degrees
respectively
$2,4,6,8$, $2,4,6,8$, $2,4,5,6,8$, or  $2,3,4,5,6,7,8$.
It follows that $H$ has, correspondingly, degrees
$2,2,4,4$, $2,2,4,4$, $2,2,4,4,5$, or $2,2,3,4,4,5,7$.
But there is no semisimple group $H$ with the last two sets of degrees:
one simple factor of $H$ would have to have maximal degree equal to an odd
number $d$ (5 or 7), hence would
be isomorphic to $SU(d)$, and hence would have all degrees
from 2 to $d$, which is not the case here. Therefore $G$ is isomorphic
to $Sp(8)$ or $Spin(9)$, and $H$ has degrees $2,2,4,4$. Since
$H$ is simply connected, Table \ref{degrees} shows that
$H$ is isomorphic to $Sp(4)^2$.

We compute that there is no free action of $Sp(4)^2$ on $Spin(9)$, and that
the only free action of $Sp(4)^2$ on $Sp(8)$ is the one-sided action given
by the standard inclusion $Sp(4)^2\subset Sp(8)$. So $M$ is the homogeneous
space $Sp(8)/Sp(4)^2$. This homogeneous space has the rational homotopy type
of $\HH\P^4\# \Ca\P^2$, but not the homotopy type, because we compute
that $w_4(Sp(8)/Sp(4)^2)=0$ whereas $w_4(\HH\P^4\# \Ca\P^2)\neq 0$, since
$w_4(\HH\P^4)\neq 0$. To compute these Stiefel-Whitney classes of the
homogeneous spaces $Sp(8)/Sp(4)^2$ and $\HH\P^4$, we can use
Singhof's approach, which works more generally for biquotients
\cite{Singhof}. Namely, say for $X=Sp(8)/Sp(4)^2$, the tangent bundle
is
$$TX=sp(8)-sp(4)_1-sp(4)_2$$
in the Grothendieck group of real vector bundles on $X$. The groups
$H^i(BSp(2n),\F_2)$ are zero for $0<i<4$, and so
$$w_4X=w_4(sp(8))-w_4(sp(4)_1)-w_4(sp(4)_2).$$
This is zero because $w_4(sp(4)_1)+w_4(sp(4)_2)$ is clearly some
$\F_2$-multiple of the sum of the generators of $H^4(BSp(4)_1,\F_2)$
and of $H^4(BSp(4)_2,\F_2)$, while
the generator of $H^4(BSp(8),\F_2)$ pulls back to the sum of the
generators of $H^4(BSp(4)_1,\F_2)$ and $H^4(BSp(4)_2,\F_2)$ and also
pulls back to zero in $H^4(X,\F_2)$.

Therefore $G$ is not simple. It must have two simple factors,
each of which contributes exactly one degree to $M$. We know that
$H$ has one more simple factor than $G$, so $H$ has three simple
factors.
We can
write $G=G_1\times G_2$ where $G_1$ contributes degree 6 to $M$
and $G_2$ contributes degree 8 to $M$, and nothing else.
By Theorems \ref{class} and \ref{top} applied to $G_1$,
there is a simple factor $H_1$ of $H$ such that either $G_1$ is
the exceptional group $G_2$, $H_1$ is $SU(2)$, and $H_1$ acts nontrivially
on $G_1$, or $H_1$ acts nontrivially
on exactly one side of $G_1$ and $G_1/H_1$ is isomorphic to one of
the homogeneous spaces  $Spin(7)/Spin(6)=G_2/SU(3)=S^6$,
$Spin(12)/Spin(11)=SU(6)/SU(5)=Sp(6)/Sp(4)=S^{11}$,
$Spin(7)/Spin(5)=UT(S^6)$,
$Spin(8)/Spin(6)$, or $Spin(8)/Spin(5)$.
 Likewise, by Theorems \ref{class} and \ref{top} applied to the
second factor $G_2$, there is a simple factor $H_2$ of $H$
which acts nontrivially on exactly one side of $G_2$ with
$G_2/H_2$ isomorphic to one of the homogeneous spaces
$Spin(9)/Spin(8)=S^8$, 
$Spin(16)/Spin(15)=SU(8)/SU(7)=Sp(8)/Sp(6)=Spin(9)/Spin(7)=S^{15}$, 
$Spin(9)/Spin(7)=UT(S^8)$,
$Spin(10)/Spin(8)$, $Spin(10)/Spin(7)$,
$Spin(9)/G_2$, or $Spin(10)/G_2$.

We see that in all the above cases $G_1/H_1$ and $G_2/H_2$, $H_1$ is never
isomorphic to $H_2$, so in particular $H_1$ is not the same simple factor of
$H$
as $H_2$. So, given $G_1/H_1$ and $G_2/H_2$, the known degrees of $M$
determine the degrees of the remaining simple factor of $H$, $H_3$.
Let us say that a given set of factors of $G$ and $H$, for example $G_1$ and
$H_1$,
{\it adds }$a$ degrees $d$ to $M$ if the integer $a$ is the number of
degrees $d$ in the given factors of $G$ minus the number of degrees $d$ in
the
given factors of $H$. From the above list, $G_1/H_1$ always adds exactly one
degree 6 to $M$, while $G_2/H_2$ adds none; so, by the known degrees of $M$,
$H_3$ has no degree 6. The only simple group with two degrees 4 is
$Spin(8)$,
which also has a degree 6; so $H_3$ has at most one degree 4.

Also, from the above list, $G_2/H_2$ adds no degree 3 to $M$, so the cases
where $G_1/H_1$ adds $-1$ degree 3 to $M$ cannot occur (otherwise $H_3$
would have $-1$ degrees 3). Thus $G_1/H_1$ is not $Spin(7)/Spin(6)=S^6$,
$G_2/SU(3)=S^6$, or $Spin(8)/Spin(6)=UT(S^7)$. Thus $G_1/H_1$
and $G_2/H_2$ both add zero degrees 3 to $M$, and so $H_3$ has no degree 3.
Next, we observe that $G_1/H_1$ adds $\geq 0$ degrees 4 to $M$, so the cases
where $G_2/H_2$ adds one degree 4 to $M$ cannot occur (otherwise $H_3$
would have at least two degrees 4). Thus $G_2/H_2$ is not
$Spin(9)/G_2$ or $Spin(10)/G_2$. Next, $G_1/H_1$ adds zero degrees 5 to $M$,
so in the cases where $G_2/H_2$ adds degree 5 to $M$, $H_3$ must have a
degree
5. Since $H_3$ has no degree 6, it follows that $H_3$ is isomorphic to
$SU(5)$, contradicting the fact that $H_3$ has no degree 3. Thus $G_2/H_2$
is not $Spin(10)/Spin(8)$ or $Spin(10)/Spin(7)$.

Finally, if $G_1/H_1$ is $Spin(8)/Spin(5)=Spin(8)/Sp(4)$, then $G_1/H_1$
adds one degree 4
 to $M$, and so, since $H_3$ has at most one degree 4, $G_2/H_2$
must add $-1$ degrees 4 to $M$; that is, $G_2/H_2$ is $Spin(9)/Spin(8)=S^8$.
Then
$H_3$ has degrees $2,4$ and so $H_3$ is isomorphic to $Sp(4)$. Here
$M$ is a biquotient $(Spin(8)\times Spin(9))/(Sp(4)\times Spin(8)\times
Sp(4))$.
Since $H_2=Spin(8)$ has finite centralizer in $G_2=Spin(9)$, the two
$Sp(4)$ factors of $H$ can only act on $Spin(9)$ on the other side from
$Spin(8)$; so $Sp(4)^2$ acts on $S^8$. Looking at the
low-dimensional
orthogonal representations of $Sp(4)^2$ shows that this action
must have a fixed point. Since $H$ acts freely on $G$,
it follows that there is a subgroup of $H$
isomorphic to $Sp(4)^2$ which acts freely on $Spin(8)$.
We compute, however, that there is no free action of $Sp(4)^2$ on $Spin(8)$.
This contradiction shows that $G_1/H_1$ is not
$Spin(8)/Spin(5)=Spin(8)/Sp(4)$.

The only remaining possibilities for $G_1/H_1$ are those
which add degree 6 to $M$ and nothing else: a homogeneous
space $Spin(12)/Spin(11)=SU(6)/SU(5)
=Sp(6)/Sp(4)=S^{11}$ or $Spin(7)/Spin(5)=UT(S^6)$ or $(G_1,H_1)=(G_2,A_1)$
(the exceptional group $G_2$).
Also, $G_2/H_2$ is either $Spin(9)/Spin(8)=S^8$, which adds
one degree 8 and subtracts one degree 4 from $M$, or else a homogeneous
space
which adds degree 8 to $M$ and nothing else: $Spin(16)/Spin(15)=SU(8)/SU(7)
=Sp(8)/Sp(6)=Spin(9)/Spin(7)=S^{15}$ or $Spin(9)/Spin(7)=UT(S^8)$.
From the known degrees of $M$, it follows that $H_3$ has only degree $2$
if $G_2/H_2$ is $S^8$ and has degrees $2,4$ otherwise. Therefore $H_3$
is isomorphic to $SU(2)$ if $G_2/H_2$ is $S^8$ and to $Sp(4)$ otherwise.

Suppose that $H_3$ is isomorphic to $Sp(4)$; we will derive a contradiction.
First, we can easily exclude the possibility that
the first factor $G_1$ is the exceptional group $G_2$. The point is that,
in this case,
neither $H_3=Sp(4)$ nor $H_2$ (from the list, above) has a nontrivial
homomorphism to the first factor $G_1$. Since $H_2\times H_3\subset H$ must
act
freely on $G_1\times G_2$, it follows that $H_2\times H_3$ acts
freely on the second factor $G_2$. This is impossible because
the list of possible spaces $G_2/H_2$ shows that $H_2$ has rank 1 less than
the second factor $G_2$,
and so $H_2\times H_3=H_2\times Sp(4)$ has rank 1 greater than the second
factor
$G_2$. The impossibility here follows from the fact
that  a finite-dimensional elliptic space $X$ has the alternating sum of its
rational homotopy groups $\chi_{\pi}(X)\leq 0$, by Halperin
(\cite{FHTbook}, p.~434). So the first factor $G_1$
is not the exceptional group $G_2$.

We continue to assume that $H_3$ is $Sp(4)$. We know that
$G_1/H_1$ is
a homogeneous
space $Spin(12)/Spin(11)=SU(6)/SU(5)
=Sp(6)/Sp(4)=S^{11}$ or $Spin(7)/Spin(5)=UT(S^6)$, and
$G_2/H_2$ is a homogeneous space $Spin(16)/Spin(15)=SU(8)/SU(7)
=Sp(8)/Sp(6)=Spin(9)/Spin(7)=S^{15}$ or $Spin(9)/Spin(7)=UT(S^8)$.
For each
possible
$G_1/H_1$ and $G_2/H_2$, we check immediately that
either $H_1$ must act trivially on $G_2$
or $H_2$ must act trivially on $G_1$. Here we use in particular that
a simple factor of $H$ must act trivially on a simple factor of $G$
isomorphic
to it, by Lemma \ref{pi3}. It follows that the biquotient
$(G_1\times G_2)/(H_1\times H_2)$ is either a $G_1/H_1$-bundle over
$G_2/H_2$ or a $G_2/H_2$-bundle over $G_1/H_1$.  Now $G_1/H_1$
is either $S^{11}$ or $UT(S^6)$, so it has the 3-local homotopy type of
$S^{11}$; and likewise $G_2/H_2$ is either $S^{15}$ or $UT(S^8)$,
so it has the 3-local homotopy type of $S^{15}$.  Therefore the biquotient
$(G_1\times G_2)/(H_1\times H_2)$ is 3-locally 11-connected. It follows
that the natural map
$$M=(G_1\times G_2)/(H_1\times H_2\times Sp(4))\arrow BSp(4)$$
is 3-locally 11-connected. Up to this point, any odd prime number would
serve in place of 3, but we now derive a contradiction by a 3-local
calculation.
Since $M$ is homotopy equivalent to $\HH\P^4\# \Ca\P^2$, the map
$\pi_8M\otimes H^8(M,\Z)\arrow \Z$ is surjective. On the other hand,
the map $\pi_8BSp(4)\times H^8(BSp(4),\Z)\arrow \Z$ has image $6\Z$.
Indeed, $H^8(BSp(4),\Z)$ is generated by $c_2^2$ and $c_4$, and $c_2^2$
is clearly zero for any $Sp(4)$-bundle over $S^8$, while $c_4$ of such a
bundle
is a multiple of $(4-1)!=6$ by Bott periodicity. Since the above map
$M\arrow BSp(4)$ is 3-locally 11-connected, we have a contradiction.

That shows that $H_3$ is not isomorphic to $Sp(4)$. Therefore $H_3$ is
isomorphic
to $SU(2)$. In this case, we know that $G_2/H_2$ is $Spin(9)/Spin(8)=S^8$.
Also, $G_1/H_1$ is either a homogeneous space $Spin(12)/Spin(11)=
SU(6)/SU(5)=Sp(6)/Sp(4)=S^{11}$ or $Spin(7)/Spin(5)=UT(S^6)$,
or else $G_1$ is the exceptional group $G_2$ and $H_1$ is $SU(2)$ acting
nontrivially on $G_1$, perhaps on both sides.

Since $M=G/H$ is homotopy equivalent to $\HH \P^4\# \Ca\P^2$, it
has the same Stiefel-Whitney classes as $\HH\P^4\# \Ca\P^2$,
by Wu (\cite{MS}, Theorem 11.14). In particular,
the sphere $S^8$ in $\HH\P^4\# \Ca\P^2$ has self-intersection number 1,
and so $w_8M$ is not in the subgroup $H^4(M,\F_2)^2$ of $H^8(M,\F_2)$.
The biquotient $N:=(G_1\times G_2)/(H_1\times H_2)$ is a principal
$SU(2)$-bundle over $M$, so its stable tangent bundle is the pullback of
that
of $M$. It follows that $w_8N\in H^8(N,\F_2)$ is not zero.

Suppose that $G_1$ is the exceptional group $G_2$. Then
$M$ is a biquotient of the form
$(G_2\times Spin(9))/(SU(2)\times Spin(8)\times SU(2))$. Here $Spin(8)$
must act trivially on $G_2$, and so we can write
$$M=(G_2\times S^8)/SU(2)^2.$$
Any torus acting by isometries on an even-dimensional sphere has a fixed
point,
since every orthogonal representation of a torus is a sum of 2-dimensional
representations and trivial representations. In particular, a maximal torus
in $SU(2)^2$ has a fixed point on $S^8$. Since $SU(2)^2$ acts freely
on $G_2\times S^8$, this torus acts freely on $G_2$. Since every element
of $SU(2)^2$ is conjugate to an element of the torus, $SU(2)^2$ acts
freely on $G_2$. Thus $M$ is a fiber bundle
$$S^8\arrow M\arrow G_2/SU(2)^2.$$
Since $S^8$ has signature zero, it follows that $M$ has signature zero,
contradicting the fact that $M$ is homotopy equivalent to
$\HH\P^4\# \Ca\P^2$. Thus the first factor $G_1$
 is not the exceptional group $G_2$.

So $G_1/H_1$ is either $Spin(12)/Spin(11)=SU(6)/SU(5)=Sp(6)/Sp(4)=S^{11}$
or $Spin(7)/Spin(5)=UT(S^6)$, and we know that
$G_2/H_2$ is $Spin(9)/Spin(8)=S^8$.
But since $w_8N$ is not zero, $N$ is not an $S^8$-bundle over
$S^{11}$ or over $UT(S^6)$. So
$H_2=Spin(8)$ must act nontrivially on
$G_1$. It follows that $G_1/H_1$ must be $Spin(12)/Spin(11)=S^{11}$.
Thus $M$ is a biquotient of the form
$$(Spin(12)\times Spin(9))/(Spin(11)\times Spin(8)\times SU(2)).$$

Here $Spin(11)$ automatically acts trivially on $Spin(9)$. So the biquotient
$$N:=(Spin(12)\times Spin(9))/(Spin(11)\times Spin(8))$$
is completely determined by the homomorphism $Spin(8)\arrow Spin(12)$,
which we know is nontrivial. (Since $Spin(11)$ has finite centralizer
on $Spin(12)$, the action of $Spin(8)$ on $Spin(12)$ must be on the other
side of $Spin(12)$ from $Spin(11)$.) There are 3 conjugacy classes of
nontrivial homomorphisms $Spin(8)\arrow Spin(12)$, each of the form
$W\oplus \R^4$, where the 8-dimensional real representation $W$
of $Spin(8)$ is either the standard representation $V$
or one of the two spin representations. Thus $N$ is the $S^{11}$-bundle
$S(W\oplus \R^4)$ over $Spin(9)/Spin(8)=S^8$. If $W$ is the standard
representation $V$ of $Spin(8)$, then the associated rank-8 vector bundle
over $S^8$ has $w_8=0$, and so the manifold $N$ would have $w_8N=0$,
contradicting what we know. So $W$ must be one of the two spin
representations
of $Spin(8)$. If necessary, we can apply an order-2 automorphism to
$Spin(8)$,
not changing the conjugacy class of the standard inclusion
$Spin(8)\arrow Spin(9)$, to arrange that $W=S^{-}$.

Thus $N$ is the $S^{11}$-bundle $S(S^{-}\oplus\R^4)$ over $S^8$. The action
of $SU(2)$ on $N$ is given by a homomorphism $SU(2)\arrow Spin(9)$
together with a homomorphism $SU(2)\arrow Spin(4)$, because $Spin(4)$
is the identity component of the centralizer of $Spin(8)$ in $Spin(12)$.
The homomorphism $SU(2)\arrow Spin(9)$ determines the action of $SU(2)$
on $S^8$ and on the $Spin(9)$-equivariant vector bundle $S^{-}$ over
$S^8$, while the homomorphism $SU(2)\arrow Spin(4)$ determines the action
of $SU(2)$ on $\R^4$. Since $SU(2)$ acts freely on $N=S(S^{-}\oplus \R^4)$,
it must act freely on $S(S^{-})=S^{15}$ and on $S(\R^4)=S^3\times S^8$.

As in the proof that $\C\P^8\#\Ca\P^2$ is not homotopy equivalent
to a biquotient, we compute that there is only one free $SU(2)$-action
on $S(S^{-}\oplus \R^4)$: $SU(2)$ must act by $SU(2)\cong Spin(3)
\subset Spin(9)$ on $S^8$ and on the bundle $S^{-}$ over $S^8$,
and by the standard representation $V_{\R}$ on $\R^4$.
Thus we have described the manifold
$M$ uniquely up to diffeomorphism. But we showed in section \ref{positive}
that exactly this manifold is diffeomorphic to $\HH\P^4\# -\Ca\P^2$.
In particular, it is not homotopy equivalent to $\HH\P^4\# \Ca\P^2$.
This completes the proof that $\HH\P^4\# \Ca\P^2$ is not homotopy
equivalent to a biquotient.

\section[The Cheeger manifold $\C\P^{4e}\# \HH\P^{2e}$]{The Cheeger manifold
$\C\P^{4e}\# \HH\P^{2e}$ is not homotopy equivalent to a biquotient}
\label{negative}

Finally, we will show that $\C\P^{4e}\# \HH\P^{2e}$, $e\geq 1$, is not
homotopy
equivalent to a biquotient $M=G/K$, completing the proof
of Theorem \ref{cheeger}.

The cohomology ring of $M$ is a complete intersection ring. The
degrees of its generators and relations 
determine the rational homotopy groups
of $M$. The result is that $G$ contributes degrees 3 and $4e$ to $M$
and $K$ contributes degrees 1 and 2, where the 1 means that
the abelianization of $K$ is isomorphic to $S^1$. Let $H$ be the commutator
subgroup of $K$, and let $N$ be the biquotient $G/H$; then $M=N/S^1$.
Here $G$ contributes degrees 3 and $4e$ to $N$, while $H$ contributes
degree 2. From the degree 2, it follows that $H$ has one more
simple factor than $G$ has.
Also, by Theorem \ref{class}, each simple factor of $G$ contributes
at least one degree to $N$. So $G$ has at most two simple factors.

Furthermore, since $M$ and $\C\P^{4e}\# \HH\P^{2e}$ are homotopy
equivalent, they have the same Stiefel-Whitney classes, by Wu (\cite{MS},
Theorem 11.14).
In particular, $w_4M$ is not in the subgroup $H^2(M,\F_2)^2$ of
$H^4(M,\F_2)$. Since $N$ is an $S^1$-bundle over $M$, the stable
tangent bundle of $N$ is the pullback of that of $M$, and so
$w_4N$ is not zero.

Suppose first that $G$ is simple. Since $H$ has one more simple factor
than $G$, $H$ has two simple factors.
We know that $G$ contributes
degrees 3 and $4e$ to $N$, and nothing else. In Theorem \ref{class},
cases (2), (3), (4) are incompatible with these degrees, and so
case (1) must hold, that is, $G$ contributes its maximal degree to $N$.
Thus $G$ has maximal degree $4e$. The only simple groups which have 3
as a degree
are the groups $SU(n)$, $n\geq 3$, and so $G$ is isomorphic to $SU(4e)$.
Since $G$ has exactly one degree $4e$, the known degrees of $N$ imply
that all simple factors of $H$ have maximal degree less than $4e$.

Suppose that $e\geq 2$. Then the second-largest degree, $4e-1$, of
$G=SU(4e)$ is killed by $H$. Let $H_1$ be a simple factor of $H$
which kills the degree $4e-1$ of $G$. We know that $H_1$ has
maximal degree less than $4e$, and so it has maximal degree $4e-1$,
which implies that $H_1$ is isomorphic to $SU(4e-1)$. Any nontrivial
homomorphism $SU(4e-1)\arrow SU(4e)$ has Dynkin index 1, and so $H_1$
must act nontrivially on only one side of $G$; otherwise
$SU(4e)$ would contribute degree 2 to $N$. So $G/H_1$
is isomorphic to the homogeneous space $SU(4e)/SU(4e-1)=S^{8e-1}$.
But then $H_1$ kills the degree 3 of $G$, which should appear in $M$.
Thus $e\geq 2$ leads to a contradiction, for $G$ simple.

This leaves the case $e=1$, with $G$ simple. As shown above, $G$
is isomorphic to $SU(4)$. Since $G$ has degrees $2,3,4$, the known
degrees of $N$ imply that $H$ has degrees $2,2$. So $H$ is isomorphic
to $SU(2)^2$. Thus $N$ is a biquotient $SU(4)/SU(2)^2$. Indeed, there is
at least one such biquotient, the homogeneous space
$UT(S^5)=Spin(6)/Spin(4)$,
which admits a free $S^1$-action. In that case, the quotient is
the 4-dimensional complex
quadric $Q^4_{\C}$, which has the rational homotopy type of
$\C\P^4\# \HH\P^2$. But no biquotient $N=SU(4)/SU(2)^2$ can have
the homotopy type of an $S^1$-bundle over $\C\P^4\# \HH\P^2$,
by the following argument.
By Singhof's
description of the tangent bundle of a biquotient \cite{Singhof},
as used in the previous section,
we have
$$w_4N=w_4(su(4))-w_4(su(2)_1)-w_4(su(2)_2).$$
We then make the convenient calculation that
$w_4(su(n))=0\in H^4(BSU(n),\F_2)$ for all even $n$. So $w_4N=0$,
which contradicts what we know about $N$. Thus we have a contradiction
from the assumption that $G$ is simple.

Thus $G$ has two simple factors. We return to the general case, $e\geq 1$.
Since $H$ has one more
simple factor than $G$, $H$ has three simple factors.
Each simple factor of $G$ contributes exactly one degree to $N$.
We can write $G=G_1\times G_2$ such that $G_1$ contributes degree 3
to $N$ and $G_2$ contributes degree $4e$, and nothing more. By
Theorems \ref{class} and \ref{top} applied to $G_1$, there
is a simple factor $H_1$ of $H$ such that either $G_1/H_1$
is isomorphic to the homogeneous space $SU(4)/Sp(4)=S^5$
(which is the same as $Spin(6)/Spin(5)$) or
$(G_1,H_1)$ is $(SU(3),SU(2))$ for some nontrivial action of $SU(2)$
on $SU(3)$.

 If $(G_1,H_1)$ is $(SU(3),SU(2))$ and $SU(2)$ acts nontrivially
on both sides of $SU(3)$, then the centralizers of both homomorphisms
$SU(2)\arrow SU(3)$  are finite or finite by $S^1$, so no other simple
factor
of $H$ can act on $G_1$. Moreover, the two nontrivial homomorphisms
$SU(2)\arrow SU(3)$ have Dynkin indices 1 and 4, by Lemma \ref{second}.
So the map $\pi_3H_1\arrow \pi_3G_1$, $\Z\arrow \Z$, is either zero or
multiplication by $1-4=-3$ or $4-1=3$. In any case, it is not surjective,
which contradicts the fact that $\pi_3H\arrow \pi_3G_1$  must be surjective,
since $\pi_3M=0$. So any $SU(2)$ factor of $H$
acts nontrivially on at most on one side of $SU(3)$. Moreover, $SU(2)$
factors
of $H$ are the only ones that can act on $G_1=SU(3)$ ($SU(3)$ factors
of $H$ are excluded by Lemma \ref{pi3}). At least one $SU(2)$ factor
must act with Dynkin index 1 rather than 4, again using that $\pi_3H\arrow
\pi_3G_1$ is surjective. Thus, if $G_1$ is $SU(3)$, then we can choose
the simple factor $H_1$ isomorphic to $SU(2)$ such that
$G_1/H_1$
is isomorphic to the homogeneous space $SU(3)/SU(2)=S^5$.

Thus, $G_1/H_1$ is either $Spin(6)/Spin(5)=S^5$ or $SU(3)/SU(2)=S^5$.
In both cases, $G_1/H_1$ adds degree 3 to $N$ and nothing else.
Next, we can apply Theorems \ref{class} and \ref{top} to the second
factor $G_2$, giving a simple factor $H_2$ with certain properties. From
the known degrees of $N=(G_1\times G_2)/(H_1\times H_2\times H_3)$,
since $G_1/H_1$ adds degree 3 to $N$ and nothing else,
$G_2/H_2$ cannot have
any degrees with multiplicity $<0$ (that is, any degrees which
occur in $H_2$ more than in $G_2$).
Given this, Theorems \ref{class} and \ref{top} imply that
either $e=1$ and $G_2$ is isomorphic to $Sp(4)$, or
there is a simple factor $H_2$ of $H$ such that
$G_2/H_2$ is isomorphic to a homogeneous space
$Spin(8e)/Spin(8e-1)=SU(4e)/SU(4e-1)=Sp(4e)/Sp(4e-2)=S^{8e-1}$,
$Spin(7)/G_2=S^7$ with $e=1$,
$Spin(9)/Spin(7)=S^{15}$ with $e=2$,
$Spin(4e+1)/Spin(4e-1)=UT(S^{4e})$
with $e\geq 2$, 
$Spin(8)/G_2=S^7\times S^7$ with $e=1$,
$Spin(9)/G_2$ with $e=2$, 
$Spin(4e+2)/Spin(4e-1)$ with $e\geq 2$,  $Spin(10)/Spin(7)$ via
the spin representation with $e=2$,
or $Spin(10)/G_2$ with $e=2$.
In these homogeneous spaces, it is clear that the simple factor $H_2$
of $H$ is different from $H_1$, which is $SU(2)$ or $Spin(5)=Sp(4)$.

When $G_2/H_2$ is one of the homogeneous spaces diffeomorphic to
$S^{8e-1}$ or $UT(S^{4e})$, 
then $G_2/H_2$
adds only one degree $4e$ to $N$ and nothing else. In those cases,
the third factor
$H_3$ of $H$ has degree 2 only, so $H_3$ is isomorphic to $SU(2)$.
When $G_2$ is isomorphic to $Sp(4)$, the product of the
two simple factors of $H$
besides $H_1$, $H_2\times H_3$,
has degrees $2,2$, and so $H_2$ and $H_3$ are both isomorphic to $SU(2)$.
When $G_2/H_2$ is $Spin(8)/G_2=S^7\times S^7$ with $e=1$ or
$Spin(9)/G_2$ with $e=2$, then $H_3$ has degrees $2,4$, and so
$H_3$ is isomorphic to $Sp(4)$. Finally, the cases where $G_2/H_2$ is
$Spin(4e+2)/Spin(4e-1)$ with $e\geq 2$, $Spin(10)/Spin(7)$ with $e=2$,
or $Spin(10)/G_2$ with $e=2$ cannot occur. In these cases, part (4)
of Theorem \ref{top} implies that $H$ must have a factor $SU(2e+1)$,
which must be the third factor $H_3$. But then $H$ has one more degree 3
than $G$, contradicting the known degrees of $N=G/H$.

The cases where $H_3$ is $Sp(4)$ are easier to analyze, so we consider them
first. Either $G_2/H_2$ is $Spin(8)/G_2=S^7\times S^7$ and $e=1$,
or $G_2/H_2$ is $Spin(9)/G_2$ and $e=2$. We know that
$G_1/H_1$ is either $SU(3)/SU(2)=S^5$ or $Spin(6)/Spin(5)=S^5$.
By the known low-dimensional representations of $H_3=Sp(4)$, the action
of $H_3$ on $G_1/H_1=S^5$ must have a fixed point. So there is a subgroup
of $H_1\times H_3$ which projects isomorphically to $H_3=Sp(4)$ and which
fixes a point in $G_1$. Also, the exceptional group $G_2$ must act
trivially on $G_1$ (which is $Spin(6)$ or $SU(3)$), so we have a subgroup
of $H$ isomorphic to $Sp(4)\times G_2$ which fixes a point in $G_1$.
Since $H$ acts freely on $G$, we have a free action of $Sp(4)\times G_2$
on the second factor $Spin(8)$ or $Spin(9)$. We compute, however,
that there is no such free action.

So we must have $H_3=SU(2)$. We know that $G_1/H_1$ is
either $SU(3)/SU(2)=S^5$ or $Spin(6)/Spin(5)=S^5$. Also, either $e=1$,
$G_2$ is isomorphic to $Sp(4)$ and $H_2$ is isomorphic to $SU(2)$,
or $G_2/H_2$ is a homogeneous space
$Spin(8e)/Spin(8e-1)= SU(4e)/SU(4e-1)=Sp(4e)/Sp(4e-2)=S^{8e-1}$,
$Spin(7)/G_2=S^7$ with $e=1$,
$Spin(9)/Spin(7)=S^{15}$ with $e=2$,
or $Spin(4e+1)/Spin(4e-1)=UT(S^{4e})$
with $e\geq 2$.

We begin with the case where $G_1/H_1$ is $SU(3)/SU(2)=S^5$ and
$G_2/H_2$ is $Spin(8e)/Spin(8e-1)=S^{8e-1}$.
This turns out to be the main
step of the whole proof; most other cases will reduce
to this one. The group
$H_2=Spin(8e-1)$ must act trivially on $G_1=SU(3)$, and so $(G_1\times G_2)/
(H_1\times H_2)$ is an $S^{8e-1}$-bundle over $S^5$. The manifold $M$
is the quotient of this bundle by a free action of  a group
$D:=K/(H_1\times H_2)$ which is
isogenous to $S^1\times SU(2)$. Since $\pi_2M$ is isomorphic to $\Z$,
$\pi_1D$ is isomorphic to $\Z$. So $D$ is isomorphic to $S^1\times
SU(2)$ or to $U(2)$.

In dimensions less than $8e-1$, $M$ clearly has the homotopy type of
a homotopy quotient $S^5\quot D$, or equivalently an $S^5$-bundle
over the classifying space $BD$. Here $D$ acts on $S^5$ through the group
$U(3)$, coming from the group $G_1=SU(3)$ together with the centralizer
$S^1$ of $H_1=SU(2)$ in $G_1$. From the cohomology ring of $M$,
the Euler class of the homomorphism $D\arrow U(3)$ must be the product
of a generator of $H^2(BD,\Z)$ with an element of $H^4(BD,\Z)$ that
generates $H^4/(H^2)^2$. Let $L$ be the standard 1-dimensional complex
representation of $S^1$, $V$ be the standard 2-dimensional
representation of $SU(2)$, and $E$ the standard 2-dimensional
representation of $U(2)$. By inspecting the low-dimensional
representations of $D$, 
it follows that
the homomorphism $D\arrow U(3)$ is isomorphic to $L^{\pm 1}\oplus
L^a\otimes V$ if $D$ is $S^1\times SU(2)$, or to $(\det E)^{\pm 1}\oplus
(\det E)^{a}\otimes E$ if $D$ is $U(2)$, for some
sign and some integer $a$. In particular, the subgroup $SU(2)=H_3$ of $D$
acts on $G_1/H_1=SU(3)/SU(2)=S^5$ by the natural inclusion
$SU(2)\arrow SU(3)$, on the other side of $G_1$ from $H_1$.

It follows that the diagonal subgroup $C$ in $H_1\times H_3
=SU(2)\times SU(2)$ has a fixed point in $SU(3)$. So $C$ must act
freely on $S^{8e-1}$. It follows that $C$ acts on $S^{8e-1}$ by the
real representation $(V_{\R})^{2e}$, where $V$ is the standard
2-dimensional complex representation of $C=SU(2)$.
In particular, the associated complex
representation is a direct sum of copies of $V$.
By the Clebsch-Gordan formula,
an irreducible complex representation of $SU(2)\times SU(2)$ whose
restriction to the diagonal subgroup is a sum of copies of $V$ must
be isomorphic to $V_1\otimes \C$ or $\C\otimes V_2$, the standard
representations of the two factors. Therefore 
the action of $H_1\times H_3=SU(2)\times SU(2)$
on $S^{8e-1}$ must be given by the real representation associated
to the complex representation
$(V_1)^{\oplus j}\oplus (V_3)^{\oplus 2e-j}$ for some $0\leq j\leq 2e$.

The manifold $M$ can be written
$$M=(SU(3)\times S^{8e-1})/((SU(2)^2\times \R)/\Z).$$
Here the subgroup $\Z$ of $SU(2)^2\times \R$ is generated by
an element of the form $(\pm 1,\pm 1,1)$, which we write as
$(e(a_0),e(b_0),1)$ for $a_0,b_0\in\{0, 1/2\}$,
where $e(t):=e^{2\pi i t}$.
The group $(SU(2)^2\times \R)/\Z$ acts on $SU(3)\times S^{8e-1}$
via homomorphisms to $SU(3)^2/Z(SU(3))$ and to $SO(8e)$.

We can assume that the first factor $SU(2)$ of $(SU(2)^2\times\R)/\Z$
acts on $SU(3)$
by the standard inclusion on the left, while the second factor $SU(2)$
acts on $SU(3)$ by the standard inclusion on the right. So the two
homomorphisms from $\R$ to $SU(3)$ both map into the centralizer
of $SU(2)$ in $SU(3)$. That is, they map $t\in\R$ to the diagonal matrices
$$(e(at),e(at),e(-2at)),(e(bt),e(bt),e(-2bt))$$
for some $a,b$. Since the generator of the above subgroup
$\Z$ in $SU(2)^2\times \R$ must map into $Z(SU(3))\cong \Z/3\subset
SU(3)^2$, $a$ and $b$ must satisfy:
\begin{align*}
a+a_0&\in \frac{1}{3}\Z\\
b+b_0&\in \frac{1}{3}\Z\\
a+a_0&\equiv b+b_0\pmod{\Z}.
\end{align*}

The homomorphism $SU(2)^2\arrow SO(8e)$ is by the real representation
$(V_1)^{\oplus j}\oplus (V_2)^{\oplus 2e-j}$ for some $0\leq j\leq 2e$.
The centralizer of this homomorphism has identity component
$Sp(2j)\times Sp(4e-2j)$, whose maximal torus is conjugate to the
center $(S^1)^{2e}$ of $U(2)^{2e}$. So we can assume that $\R$ maps
into this center. Thus the homomorphism from $(SU(2)^2\times\R)/\Z$
to $SO(8e)$ is
the direct sum of $2e$ homomorphisms to $U(2)$,
of the form $(A,B,t)\mapsto e(c_it)A$ for $1\leq i\leq j$
and $(A,B,t)\mapsto e(d_it)B$ for $j+1\leq i\leq 2e$. Since this
homomorphism is trivial on the subgroup $\Z$, the numbers $c_i$ and
$d_i$ must satisfy $c_i+a_0\in\Z$ and $d_i+b_0\in\Z$.

The action of $(SU(2)^2\times \R)/\Z$ on
$SU(3)\times S^{8e-1}$ is free. We compute that this means that,
for all $1\leq i\leq j$,
\begin{align*}
a+2b-c_i=\pm 1\\
a+2b+c_i=\pm 1\\
-2a+2b=\pm 1,
\end{align*}
and for all $j+1\leq i\leq 2e$,
\begin{align*}
2a+b-d_i=\pm 1\\
2a+b+d_i=\pm 1\\
-2a+2b=\pm 1.
\end{align*}
If $j\geq 1$, then the first two equations imply that
either $a+2b=0$ and $c_i=\pm 1$, or $a+2b=\pm 1$ and $c_i=0$,
for $1\leq i\leq j$. In particular, $c_i$ is an integer,
which implies (since $c_i+a_0\in\Z$) that $a_0$ is 0, not $1/2$.
Since $a+a_0\in (1/3)\Z$, it follows in particular that $a$ is 2-integral;
since $-2a+2b=\pm 1$, $b$ is not 2-integral. Therefore $b_0$ is
$1/2$, not 0. But if $j\leq 2e-1$, then we would get the opposite
conclusion (that $a_0=1/2$ and $b_0=1$) from the second three
equations above. So we must have either $j=0$ or $j=2e$.
After switching the two $SU(2)$ factors if necessary, we can assume
that $j=2e$. That is, $SU(2)^2$ acts on $S^{8e-1}$ by the
real representation $(V_1)^{2e}$. Also, we have $a_0=0$ and $b=1/2$,
which means that the subgroup $\Z$ of $SU(2)^2\times \R$
is generated by $(1,-1,1\in\R)$.

Since the second factor $SU(2)$ acts only on $SU(3)$, we can rewrite $M$
as
$$M=(S^5\times S^{8e-1})/(SU(2)\times S^1).$$
Here we have used that the quotient of $(SU(2)^2\times\R)/\Z$ by the
second copy of $SU(2)$ is isomorphic to $SU(2)\times S^1$. The action
of $SU(2)\times S^1$ on $S^{8e-1}$ is given by the complex representation
$\oplus _{i=1}^{2e} V\otimes L^{c_i}$. Also, we compute that the action
of $SU(2)\times S^1$ on $S^5$ is by the complex representation
$(V\otimes L^{a+2b})\oplus L^{-2a+2b}$. As computed above, 
$-2a+2b=\pm 1$. And either
$a+2b=\pm 1$ and $c_i=0$ for all $i$, or $a+2b=0$ and $c_i=\pm 1$
for all $i$. 

In the first case, where $c_i=0$,
we can conjugate the action of $SU(2)\times S^1$ on $S^5$
in the orthogonal group $O(6)$ to make $a+2b$
and $-2a+2b$ equal to 1,
rather than $-1$. Thus $M$ is the manifold
$$(S((V\otimes L)\oplus L)\times S(V^{\oplus 2e}))/(SU(2)\times S^1).$$
So $M$ is a $\C\P^2$-bundle over
$\HH\P^{2e-1}$, and so $M$ has signature zero, contradicting that
$M$ is homotopy equivalent to $\C\P^{4e}\# \HH\P^{2e}$. In fact,
this $\C\P^2$-bundle over $\HH\P^{2e-1}$ is the one diffeomorphic to
$\C\P^{4e}\# -\HH\P^{2e}$, as mentioned in section \ref{positive}.
In the second case, where $c_i=\pm 1$, we can conjugate the homomorphisms
from $SU(2)\times S^1$ to the orthogonal groups
$O(6)$ and $O(8e)$ to make $-2a+2b=1$ and
$c_i=1$ for all $i$.
Thus, $M$ is the manifold
$$(S(V\oplus L)\times S((V\otimes L)^{\oplus 2e}))/(SU(2)\times S^1).$$
Again, we showed in section \ref{positive}
 that this manifold is diffeomorphic to $\C\P^{4e}\# -\HH\P^{2e}$.
Thus it is not homotopy equivalent to $\C\P^{4e}\# \HH\P^{2e}$.
This completes
the proof that we cannot have $G_1/H_1$ equal to $SU(3)/SU(2)=S^5$
and $G_2/H_2$ equal to $Spin(8e)/Spin(8e-1)=S^{8e-1}$.

We next consider the case where $G_1/H_1$ is $Spin(6)/Spin(5)=S^5$ and
$G_2/H_2$ is $Spin(8e)/Spin(8e-1)=S^{8e-1}$.
A first observation is that $Spin(5)$ has finite centralizer
in $Spin(6)$, and so all factors of $K$ except $Spin(5)$ act on the other
side
of $Spin(6)$ from $Spin(5)$. Likewise, all factors of $K$ except
$Spin(8e-1)$ act on the other side of $Spin(8e)$ from $Spin(8e-1)$.
Furthermore, $Spin(8e-1)$ must act trivially on $Spin(6)$,
and so $(G_1\times G_2)/(H_1\times H_2)$ is an $S^{8e-1}$-bundle over
$S^5$, which we can write as $(Spin(6)\times S^{8e-1})/Spin(5)$.

The homomorphism $K\arrow (G\times G)/Z(G)$ which defines
the action of $K$ on $G$ gives a homomorphism
from $D:=K/(H_1\times H_2)$ to $SO(6)$. Here $D$ is isogenous
to $SU(2)\times S^1$ and has fundamental group isomorphic to $\Z$, so
$D$ is isomorphic to $SU(2)\times S^1$ or to $U(2)$.
In dimensions less than $8e-1$, $M$ has the homotopy
type of the homotopy quotient $S^5\quot D$ defined by the homomorphism
$D\arrow SO(6)$, or equivalently of an $S^5$-bundle over the classifying
space $BD$. 
From the cohomology ring of $M$, the Euler class of the homomorphism
$D\arrow SO(6)$
must be the product of some generator of $H^2(BD,\Z)$ with
some element of $H^4(BD,\Z)$ which generates $H^4/(H^2)^2$.
By inspecting the low-dimensional real representations of $D$, it follows
that the homomorphism $D\arrow SO(6)$ is the real representation
associated to a 3-dimensional complex representation of $D$. For
$D=S^1\times
SU(2)$, write $L$ for the standard 1-dimensional complex representation
of $S^1$ and $V$ for the standard 2-dimensional complex
representation of $SU(2)$.
In order to have an Euler class of the form above,
the homomorphism $D\arrow SO(6)$ must come from the complex
representation $L^{\pm 1}+L^b\otimes V$ for some sign and some integer
$b$. The Euler class of this representation is $\pm xy$,
where we let $x=c_1L$ and $y=c_2V+b^2x^2$.
Likewise, for $D=U(2)$, let $E$ be the standard 2-dimensional complex
representation of $D$. In order to have an Euler class of the form
above,
the homomorphism
$D\arrow SO(6)$ must come from the complex representation
$(\det E)^{\pm 1}+(\det E)^d\otimes E$ for some sign and some integer $d$.
The Euler class in $H^6BU(2)$ of this representation is $\pm uv$,
where we let $u=c_1E$ and $v=c_2E+(d^2+d)u^2$.

Because the homomorphism $D\arrow SO(6)$ factors through $U(3)$,
we can replace $G_1/H_1=Spin(6)/Spin(5)=S^5$  by
$G_1/H_1=\widetilde{U}(3)/\widetilde{U}(2)$, in our description of $M$
as a biquotient. Here $\widetilde{U}(n)$ denotes the inverse image of
$U(n)\subset SO(2n)$ in $Spin(2n)$. We can then apply the proof of
Lemma \ref{simpconn} to replace $\widetilde{U}(3)$ by $SU(3)$.
Thus we have reduced to the case where $G_1/H_1$ is $SU(3)/SU(2)=S^5$ and
$G_2/H_2$ is $Spin(8e)/Spin(8e-1)=S^{8e-1}$. But we have shown that
the latter case cannot occur. This completes the proof that
$G_2/H_2$ cannot be $Spin(8e)/Spin(8e-1)=S^{8e-1}$, either when
$G_1/H_1$ is $Spin(6)/Spin(5)$ or when it is $SU(3)/SU(2)$.

The situation now is as follows. First, we know that $G_1/H_1$ is
either $SU(3)/SU(2)=S^5$ or $Spin(6)/Spin(5)=S^5$. Also, either
$e=1$, $G_2$ is isomorphic to $Sp(4)$ and $H_2$ is isomorphic to $SU(2)$, or
$G_2/H_2$ is a homogeneous space
$Spin(8e)/Spin(8e-1)= SU(4e)/SU(4e-1)=Sp(4e)/Sp(4e-2)=S^{8e-1}$,
$Spin(7)/G_2=S^7$ with $e=1$,
$Spin(9)/Spin(7)=S^{15}$ with $e=2$, or
$Spin(4e+1)/Spin(4e-1)=UT(S^{4e})$
with $e\geq 2$.
We have shown that $G_2/H_2$ cannot be
$Spin(8e)/Spin(8e-1)=S^{8e-1}$, either when $G_1/H_1$ is
$Spin(6)/Spin(5)$ or when it is $SU(3)/SU(2)$.
Most other cases reduce to this one.
Namely, suppose  that $G_2/H_2$ is diffeomorphic to a sphere $S^{8e-1}$
(noting that in all these cases $G_2\times Z_{G_2}(H_2)$
acts by isometries in the usual metric
on the sphere), and that $H_2$ acts
trivially on $G_1$. Then we can simply replace $G_2/H_2$ by
$Spin(8e)/Spin(8e-1)=S^{8e-1}$  without changing the biquotient
$M$. Since $G_1$ is small, namely $Spin(6)$ or $SU(3)$, the hypothesis
that $H_2$ acts trivially on $G_1$ is automatic for most of the pairs
$G_2/H_2$.
Namely this holds when there is no nontrivial homomorphism $H_2\arrow G_1$,
or also when $H_2$ is isomorphic to $G_1$ by Lemma \ref{pi3}.

The cases not covered by this argument are: $G_1/H_1$ is $SU(3)/SU(2)=S^5$
or $Spin(6)/Spin(5)=S^5$, $G_2$ is isomorphic to $Sp(4)$,
and $H_2$ is isomorphic to $SU(2)$,
where $e=1$;
$G_1/H_1$ is $Spin(6)/Spin(5)$ and $G_2/H_2$ is $SU(4)/SU(3)=S^{15}$ with $H_2$
acting nontrivially on $G_1$, where $e=2$;  or
$G_1/H_1$ is $SU(3)/SU(2)$ or $Spin(6)/Spin(5)$ and $G_2/H_2$
is $Spin(4e+1)/Spin(4e-1)=UT(S^{4e})$ with $e\geq 2$.

The last case, where $G_2/H_2$ is $UT(S^{4e})$
with $e\geq 2$, is easy to exclude.
Let $N=(G_1\times G_2)/(H_1\times H_2\times H_3)$, so that $M=N/S^1$.
Because
the groups involved in $N$ are simply connected,
$N$ is 2-connected, and so $N$ is the $S^1$-bundle over $M$ corresponding
to a generator of $H^2(M,\Z)$. Since $M$ has the integral cohomology ring of
$\C\P^{4e}\# \HH\P^{2e}$, the spectral sequence of this $S^1$-bundle
shows that $N$ has the integral cohomology ring of $S^5\times
\HH\P^{2e-1}$.
Next, let $Y$ be the manifold $(G_1\times G_2)/(H_1\times H_2)$, which
is an $SU(2)$-bundle over $N$ because $H_3=SU(2)$. Because $H_2=Spin(4e-1)$
and $e\geq 2$, $H_2$ acts trivially on $G_1$ (which is $Spin(6)$ or
$SU(3)$),
and so $Y$ is a $UT(S^{4e})$-bundle over $S^5$. In particular, $Y$
is 4-connected. Also, the spectral sequence computing the cohomology
of $Y$ collapses for degree reasons, and so $Y$ has 2-torsion in its
cohomology
because $UT(S^{4e})$ does. But $Y$ is also an $S^3$-bundle over $N$,
and because $Y$ is 4-connected, the Euler class of this bundle must
be a generator of $H^4(N,\Z)$. So the spectral sequence of this $S^3$-bundle
shows that $Y$ has the integral cohomology ring of $S^5\times S^{8e-1}$.
This contradicts the fact
that $Y$ has 2-torsion. So this case, $G_2/H_2=UT(S^{4e})$ with $e\geq 2$,
does not occur.

Next, we consider the case where $G_1/H_1$ is $Spin(6)/Spin(5)=S^5$
and $G_2/H_2$ is $SU(4)/SU(3)=S^{15}$, with $SU(3)$ acting nontrivially
on $Spin(6)$. Here $e=2$. The point is that any nontrivial action
of $SU(3)$ on $S^5$ is isomorphic to the standard action, and hence
is transitive. Thus $H$ acts transitively on the factor $G_1$ of $G$,
contrary
to Convention \ref{simp}.

This completes the proof that $\C\P^{4e}\#\HH\P^{2e}$ is not homotopy
equivalent to a biquotient for all $e\geq 2$. Ironically, the hardest case
of all is the case $e=1$, that is, the proof that $\C\P^4\# \HH\P^2$ is not
homotopy equivalent to a biquotient.

For $e=1$, it remains to consider the case where $G_1/H_1$ is either
$SU(3)/SU(2)=S^5$ or $Spin(6)/Spin(5)=S^5$, $G_2$ is isomorphic to $Sp(4)$,
and $H_2$ is isomorphic to $SU(2)$.
Thus the 9-manifold $N$ is a biquotient of the
form
$$(SU(3)\times Sp(4))/SU(2)^3$$
or
$$(Spin(6)\times Sp(4))/(Spin(5)\times SU(2)^2),$$
where the first factor $H_1$ acts by the standard inclusion on one side
of $G_1$ and trivially on the other side of $G_1$.
 Also, if $H_2$ or $H_3$ (each isomorphic to $SU(2)$)
acts trivially on $G_1$,
and acts trivially on
one side of $Sp(4)$ and by the standard inclusion $V\oplus \C^2$ on the
other,
then we can replace the quotient $Sp(4)/SU(2)=S^7$ by $Spin(8)/Spin(7)=S^7$
and thus reduce to an earlier case. So we can assume that neither $H_2$
nor $H_3$ has these properties.

We checked earlier that the Stiefel-Whitney class $w_4N$ is not zero.
It is conven\-ient to observe now that the Pontrjagin class $p_1(\HH\P^2)$
is $2z$, where $z$ is a generator of $H^4(\HH\P^2,\Z)\cong
\Z$. Furthermore, Wu showed that the first Pontrjagin class of
a closed manifold is a homotopy invariant modulo 12 \cite{Wu}.
Since $M$ is homotopy equivalent to $\C\P^4\# \HH\P^2$,
$p_1M$ in $H^4M/(H^2M)^2\cong \Z$ is 2 times the class of some generator,
modulo 12.
Therefore, the $S^1$-bundle $N$ over $M$ has
$p_1N$ equal to 2 times the class of some generator of $H^4N\cong \Z$,
modulo 12.

\begin{lemma}
\label{free}
If $H_1$ acts trivially on $G_2=Sp(4)$, then $H_2\times H_3=SU(2)^2$
does not act freely on $G_2=Sp(4)$.
\end{lemma}

{\bf Proof. }Suppose that $H_1$ acts trivially on $Sp(4)$ and that
$H_2\times H_3$ acts freely on $Sp(4)$. Since $H_1$ acts trivially
on $Sp(4)$, we can enlarge $G_1$ and $H_1$ if necessary to make
$G_1/H_1$ equal to $Spin(6)/Spin(5)$, rather than $SU(3)/SU(2)$.
The quotient
$Sp(4)/(H_2\times H_3)$
is diffeomorphic to $S^4$, by Lemma \ref{sp4}. Also, $N$ is the
$S^5$-bundle over $S^4$ associated to some homomorphism
from $H_2\times H_3$ to $Spin(6)$. By our knowledge of $p_1N$,
$p_1$ of this $S^5$-bundle in $H^4(S^4,\Z)\cong\Z$
must be 2 times some generator, modulo 12.

By Lemma \ref{sp4}, after switching $H_2$ and $H_3$ and switching
the two sides of $Sp(4)$ if necessary, $H_2\times H_3$ acts on $Sp(4)$
by either $(V_2\oplus V_3,\C^4)$ or $(V_2\oplus \C^2,(V_3)^{\oplus 2})$.
First suppose that $H_2\times H_3$ acts on $Sp(4)$ by $(V_2\oplus
V_3,\C^4)$. The conjugacy classes of homomorphisms from $H_2\times H_3$
to $Spin(6)$ have complexifications: $\C^6$, $(S^2V_2)^{\oplus 2}$,
$S^2V_2+S^2V_3$, $(S^2V_3)^{\oplus 2}$, $(V_2)^{\oplus 2}\oplus \C^2$, 
$(V_3)^{\oplus 2}\oplus \C^2$, $V_2\otimes V_3\oplus \C^2$, $S^4V_2\oplus\C$,
and $S^4V_3\oplus \C$. The given homomorphism from $H_2\times H_3$
to $Spin(6)$ must have first Pontrjagin class (that is, $-c_2$
of the complexification) equal to 2 times some
generator of $H^4(S^4,\Z)$ modulo 12, where both $c_2V_2$ and $c_2V_3$
represent the same generator of $H^4(S^4,\Z)$. This only occurs
when the complexification is $(V_2)^{\oplus 2}\oplus \C^2$ or
$(V_3)^{\oplus 2}\oplus \C^2$. But then one of  $H_2$ or $H_3$ acts
trivially on $Spin(6)$ and by $(V\oplus\C^2,\C^4)$ on $Sp(4)$,
contrary to what we arranged before the lemma.

It remains to consider the case where $H_2\times H_3$ acts on $Sp(4)$
by $(V_2\oplus\C^2,(V_3)^{\oplus 2})$. In this case, $c_2V_3$ represents
a generator of $H^4(S^4,\Z)$, while $c_2V_2$ is 2 times that generator.
Going through the list of possible homomorphisms from $H_2\times
H_3$ to $Spin(6)$ again, we find that the only one whose
Pontrjagin class in $H^4(S^4,\Z)$ is 2 times a generator, modulo 12,
is the one with complexification $(V_3)^{\oplus 2}\oplus\C^2$.
But then $H_2$ acts trivially on $Spin(6)$ and by
$(V\oplus\C^2,\C^4)$ on $Sp(4)$, contrary to what we arranged
before the lemma. \qed

We can now show that the case where $G_1/H_1$ is $Spin(6)/Spin(5)$
does not occur. In this case, we know that $H_1$ acts trivially on
$Sp(4)$ (which is isomorphic to $Spin(5)$), by Lemma \ref{pi3}.
By Lemma \ref{free}, $H_2\times H_3$ does not act freely on $Sp(4)$.
On the other hand, the action of $H_2\times H_3$ on $Spin(6)$
is given by one of the homomorphisms listed in the proof of
Lemma \ref{free}. In all these cases, we check immediately that
the maximal torus $(S^1)^2$ in $H_2\times H_3=SU(2)^2$ has a fixed
point in $Spin(6)/Spin(5)=S^5$. Therefore, there is a subgroup
of $H=H_1\times H_2\times H_3$ which projects isomorphically to
$(S^1)^2$ in $H_2\times H_3$ and which has a fixed point on
$Spin(6)$. Since $H$ acts freely on $G$, this subgroup acts
freely on $Sp(4)$. Since $H_1$ acts trivially on $Sp(4)$, this
means that the maximal torus $(S^1)^2$ in $H_2\times H_3$ acts
freely on $Sp(4)$. Since every element of $H_2\times H_3$ is conjugate
to an element of the maximal torus, it follows that
$H_2\times H_3$ acts freely on $Sp(4)$, contradicting what we have
shown. Thus the case where $G_1/H_1$ is $Spin(6)/Spin(5)$ does
not occur.

It remains to consider the case where $G_1/H_1$ is $SU(3)/SU(2)=S^5$.
That is,
the 9-manifold $N$ is a biquotient
$$(SU(3)\times Sp(4))/SU(2)^3,$$
where we know that $H_1$ acts on $SU(3)$ by the standard inclusion on one
side, and we arranged earlier that neither $H_2$ nor $H_3$
acts trivially on one side of $Sp(4)$ and by the standard inclusion on the
other.

Suppose first that $H_2$ and $H_3$ act trivially on $G_1$. Then $H_2\times
H_3
\cong SU(2)^2$
must act freely on $Sp(4)$. By Lemma \ref{sp4},
up to switching $H_2$ and $H_3$ and switching the two sides
of $Sp(4)$, $H_2\times H_3$ acts on $Sp(4)$ by the homomorphisms
$(V_2\oplus \C^2,V_3\oplus V_3)$ or $(V_2\oplus V_3,\C^4)$. Therefore,
one of the factors $H_2$ or $H_3$ acts on $Sp(4)$ trivially on one side
and by the standard inclusion on the other, as well as acting trivially on
$SU(3)$.
This contradicts what we arranged earlier.

So one of $H_2$ or $H_3$ acts nontrivially on $G_1$. Switching $H_2$ and
$H_3$ if necessary, we can assume that $H_2$ acts nontrivially on $G_1$.
Since the centralizer of $H_1=SU(2)$
in $G_1=SU(3)$ is only finite by $S^1$,
$H_2$ acts only on the other side of $G_1$ from $H_1$. It must act
on $G_1=SU(3)$ by the homomorphism $V_2\oplus \C$ or $S^2V_2$.
Thus the centralizer of $H_2$ in $G_1$ is at most finite by $S^1$.
Since the centralizer of $H_1$ on the other side of $G_1$
is finite by $S^1$, $H_3=SU(2)$ must act trivially on $G_1$.
By what we arranged earlier, it follows that
$H_3$ does not act on $Sp(4)$ trivially on one side and by the
standard inclusion on the other.

Suppose first that
$H_2$ acts on $SU(3)$ by $V_2\oplus \C$. Then $H_2$ has a fixed
point on $S^5$. More precisely, we see that the diagonal subgroup
$\Delta_{12}\cong SU(2)$ in $H_1\times H_2$ has a fixed point in $SU(3)$.
Since $H_3$ acts trivially on $G_1=SU(3)$, it follows that
$\Delta_{12}\times H_3\cong SU(2)^2$ acts freely on $G_2=Sp(4)$. By the
classification of free actions of $SU(2)^2$ on $Sp(4)$ in Lemma \ref{sp4},
together with the fact that $H_3$ does not
act trivially on one side of $Sp(4)$ and by the standard inclusion
on the other, $\Delta_{12}\times H_3$ must act on $Sp(4)$
by the homomorphism $V_{12}\oplus\C^2$ on one side and $V_3\oplus V_3$
on the other. Here $V_{12}$ denotes the standard representation
of $\Delta_{12}=SU(2)$. By the Clebsch-Gordan formula, it follows that
$H_1\times H_2\times H_3=SU(2)^3$ acts on $Sp(4)$ by
$V_1\oplus\C^2$ or $V_2\oplus \C^2$ on one side and
by $V_3\oplus V_3$ on the other. Switching $H_1$ and $H_2$ if necessary
(which we can do, since they act the same way on $G_1=SU(3)$),
we can assume that $H_1\times H_2\times H_3=SU(2)^3$ acts on $Sp(4)$ by
$V_2\oplus \C^2$ on one side and by $V_3^{\oplus 2}$ on the other.
Thus, $H_1$ acts trivially on $Sp(4)$ and $H_2\times H_3$ acts freely
on $Sp(4)$, which contradicts Lemma \ref{free}.

Therefore $H_2=SU(2)$ must act on $G_1=SU(3)$ by the homomorphism
$S^2V_2$. By Singhof's description of the tangent bundle of
a biquotient \cite{Singhof}, $w_4N$ is given by
\begin{align*}
w_4N &= w_4su(3)+w_4sp(4)-w_4su(2)_1-w_4su(2)_2-w_4su(2)_3\\
&= w_4su(3)+w_4sp(4),
\end{align*}
using that $w_4su(2)=0$ in $H^4(BSU(2),\F_2)$. Furthermore, we know
that $H_2$ acts on one side of $SU(3)$ by $S^2V_2$, and no other factor
of $H$ acts on that side of $SU(3)$. Since $c_2(S^2V_2)=4c_2V_2$,
the generator $c_2$ of $H^4(BSU(3),\F_2)$ pulls back to
$4c_2V_2=0$ in $H^4(N,\F_2)$. Therefore
$$w_4N=w_4sp(4).$$
We know that $H_3$ acts trivially on $SU(3)$, so it must act freely on
$Sp(4)$.
Furthermore, we know that it does not act trivially on one side of $Sp(4)$
and
by the standard inclusion $V_3\oplus \C^2$ on the other.
By Lemma \ref{sp4}, $H_3$ must act on at least one side of $Sp(4)$ by
the homomorphism $V_3^{\oplus 2}$ or $S^3V_3$. These two homomorphisms have
centralizers in $Sp(4)$ which are finite by $S^1$ or finite, so no other
simple factor
 of $H$ acts on the same side of $Sp(4)$. Since $c_2(V_3^{\oplus 2})=2c_2V_3$
and $c_2(S^3V_3)=10c_2V_3$, the generator $c_2$ of $H^4(BSp(4),\F_2)$
pulls back to $2c_2V_3$ or $10c_2V_3$ in $H^4(N,\F_2)$, thus to zero.
Therefore
$$w_4N=0,$$
contradicting what we know about $N$. This completes the proof that
$\C\P^4\# \HH\P^2$ is not homotopy equivalent to a biquotient.
Theorem \ref{cheeger} is proved. \qed


\small \sc DPMMS, Wilberforce Road,
Cambridge CB3 0WB, England.

b.totaro@dpmms.cam.ac.uk

\begin{thebibliography}{99}

\bibitem{AA} J.~F.~Adams and M.~Atiyah. $K$-theory and the Hopf invariant.
{\it Quart.\ J.\ Math.\ }{\bf 17 }(1966), 31--38.

\bibitem{Aguilar} R.~Aguilar. Symplectic reduction and the complex
homogeneous Monge-Amp\`ere equation.
{\it Ann.\ Global Anal.\ Geom.\ }{\bf 19 }(2001), 327--353.

\bibitem{Berger} M.~Berger. Les vari\'et\'es riemanniennes homog\`enes
normales simplement connexes \`a courbure strictement positive.
{\it Ann.\ Scuola Norm.\ Sup.\ Pisa }{\bf 15 }(1961), 179--246.

\bibitem{Borel} A.~Borel. Le plan projectif des octaves et les
sph\`eres comme espaces homog\`enes.
{\it C.~R.\ Acad.\ Sci.\ }{\bf 230 }(1950), 1378--1380.

\bibitem{Bourbaki} N.~Bourbaki. {\it Groupes et alg\`ebres de Lie,
Chapitres 4, 5 et 6. }Masson (1981).

\bibitem{Cheeger} J.~Cheeger. Some examples of manifolds of nonnegative
curvature. {\it J.~Diff.\ Geo.\ }{\bf 8 }(1973), 623--628.

\bibitem{Dold} A.~Dold. Erzeugende der Thomschen Algebra ${\mathfrak R}$.
{\it Math.\ Zeit.\ }{\bf 65 }(1956), 25-35.

\bibitem{Dynkin} E.~Dynkin. Semisimple subalgebras of semisimple
Lie algebras. {\it Mat.\ Sbornik }{\bf 30 }(1952), 349--462; AMS
Translations
(2) {\bf 6 }(1957), 111--244.

\bibitem{Eschenburgexample} J.~Eschenburg. New examples of manifolds with
strictly positive curvature. {\it Invent.\ Math.\ }{\bf 66 }(1982),
469--480.

\bibitem{EschenburgHab} J.~Eschenburg. Freie isometrische Aktionen auf
kompakten Lie-Gruppen mit positiv gekr\"ummten Orbitr\"aumen.
{\it Schriften der Math.\ Universit\"at M\"unster }{\bf 32 }(1984).

\bibitem{Eschenburgcoho} J.~Eschenburg. Cohomology of biquotients.
{\it Manu.\ Math.\ }{\bf 75 }(1992), 151--166.

\bibitem{FHTell} Y.~F\'elix, S.~Halperin, and J.-C.~Thomas. Elliptic
spaces II. {\it Ens.\ Math.\ }{\bf 39 }(1993), 25--32.

\bibitem{FHTbook} Y.~F\'elix, S.~Halperin, and J.-C.~Thomas. {\it
Rational homotopy theory. }Springer (2001).

\bibitem{GO} V.~Gorbatsevich and A.~Onishchik. Lie transformation groups,
95--226. {\it Lie groups and Lie algebras I}, ed.\ A.~Onishchik,
Encyclopaedia
of Math. Sciences, v.~20, Springer (1993).

\bibitem{GM} D.~Gromoll and W.~Meyer. An exotic sphere with
nonnegative sectional curvature. {\it Ann.\ Math.\ }{\bf 100 }(1974),
401--406.

\bibitem{GZ} K.~Grove and W.~Ziller. Curvature and symmetry of Milnor
spheres. {\it Ann.\ Math.\ }{\bf 151 }(2000), 1--36.

\bibitem{KZ} V.~Kapovitch and W.~Ziller.
Biquotients with singly generated rational cohomology.
Preprint (2002).

\bibitem{Lambrechts} P.~Lambrechts. The Betti numbers of the free loop
space of a connected sum. {\it J.\ LMS }{\bf 64 }(2001), 205--228.

\bibitem{MS} J.~Milnor and J.~Stasheff.
{\it Characteristic classes. }Princeton (1974).

\bibitem{Onishchikinclusion} A.~Onishchik. Inclusion relations between
transitive compact transformation groups.
{\it Trudy Moskov.\ Mat.\ Obshch.\ }{\bf 11 }(1962), 199--242;
{\it AMS Translations }(2) {\bf 50 }(1966), 5--58.

\bibitem{Onishchikreflection} A.~Onishchik. Remark on invariants of groups
generated by reflections. {\it Selecta Math.\ Sov.\ }{\bf 3 }(1983/84),
239--241.

\bibitem{Onishchikbook} A.~Onishchik. {\it Topology of transitive
transformation groups. }Johann Ambrosius Barth (1994).

\bibitem{Paternain} G.~Paternain. Differentiable structures with zero
entropy on simply connected 4-manifolds. {\it Bol.\ Soc.\ Brasil.\
Mat.\ }{\bf 31 }(2000), 1--8.

\bibitem{PT} A.~Petrunin and W.~Tuschmann. Diffeomorphism finiteness,
positive pinching, and second homotopy. {\it Geom.\ Funct.\
Anal.\ }{\bf 9 }(1999), 736--774.

\bibitem{Singhof} W.~Singhof. On the topology of double coset manifolds.
{\it Math.\ Ann.\ }{\bf 297 }(1993), 133--146.

\bibitem{Totarocomplex} B.~Totaro. Complexifications of nonnegatively
curved manifolds. {\it J.\ Eur.\ Math.\ Soc., }to appear.

\bibitem{Totarocurv} B.~Totaro. Curvature, diameter, and quotient manifolds.
math.DG/0209173

\bibitem{Wu} W.-T. Wu. On Pontrjagin classes III. {\it AMS
Translations }(2) {\bf 11 }(1959), 155--172.

\end{thebibliography}
\end{document}